\documentclass[12pt]{article}
\usepackage{amsmath}
\usepackage{amssymb}
\usepackage{latexsym}
\usepackage{graphicx}
\usepackage{amsfonts}
\usepackage{amscd}

\def\black{\ {\hbox{\vrule width 4pt height 4pt depth
0pt}}}

\def\fine{\ \black\vskip.4truecm}
\def\Proof{{\sl Proof.}\quad}

\def\set#1{\{\,#1\,\}}

\let\Im\undefined
\let\End\undefined

\DeclareMathOperator{\End}{End}
\DeclareMathOperator{\Ker}{Ker}

\DeclareMathOperator{\Im}{Im}
\DeclareMathOperator{\Ann}{Ann}

\def\Ext#1#2#3#4{\mathop{{\mathrm {Ext}}^{#1}_{#2}(#3,#4)}}

\def\t{\hbox{\rm t}}

\def\Z{\mathbb Z}

\def\Q{\mathbb Q}
\def\T{\mathbb T}

\newtheorem{DEF}{Definition}[section]

\newtheorem{EXS*}[DEF]{Examples}
\newtheorem{EX*}[DEF]{Example}
\newtheorem{REM*}[DEF]{Remark}

\def\dualita#1#2{\mathrel{
                 \mathop{\vcenter{
                 \offinterlineskip
                 \hbox to 1.2truecm{\rightarrowfill}
                 \hbox to 1.2truecm{\leftarrowfill}}}%
                 \limits_{#2}^{#1}}}

\begin{document}
\newtheorem{theorem}{Theorem}[section]
\newtheorem{remark}[theorem]{Remark}
\newtheorem{definition}[theorem]{Definition}
\newtheorem{notation}[theorem]{Notation}
\newtheorem{Sh-ele}[theorem]{Shelah's Elevator}
\newtheorem{fprop}[theorem]{Freeness Proposition}
\newtheorem{mtheorem}[theorem]{Main Theorem}
\newtheorem{steplemma}[theorem]{Step Lemma}
\newtheorem{mlemma}[theorem]{Main Lemma}
\newtheorem{defobs}[theorem]{Definition-Observation}
\newtheorem{blackbox}[theorem]{The General Black Box}
\newtheorem{weakdia}[theorem]{The Weak Diamond Principle }
\newtheorem{sblackbox}[theorem]{The Strong Black Box}
\newtheorem{wllemma}[theorem]{Wald-{\L}o{\'s} Lemma}
\newtheorem{pigeonlemma}[theorem]{Pigeon-hole Lemma}
\newtheorem{observation}[theorem]{Observation}
\newtheorem{proposition}[theorem]{Proposition}
\newtheorem{lemma}[theorem]{Lemma}
\newtheorem{corollary}[theorem]{Corollary}
\newtheorem{construction}[theorem]{Construction}
\newtheorem{example}[theorem]{Example}
\newtheorem{recognition}[theorem]{Recognition Lemma}
\renewcommand{\labelenumi}{(\roman{enumi})}
\newcommand{\dach}[1]{\hat{\vphantom{#1}}}
\numberwithin{equation}{section}
\def\Pf{\smallskip\goodbreak{\sl Proof. }}
\def\Fp{\vadjust{}\penalty200 \hfill
\lower.3333ex\hbox{\vbox{\hrule\hbox{\vrule\phantom{\vrule height
6.83333pt depth 1.94444pt width 8.77777pt}\vrule}\hrule}}
\ifmmode\let\next\relax\else\let\next\par\fi \next}
\def\qa{(A_\a, G_\a, \s_\a)}
\def\qb{(A_\b, G_\b, \s_\b)}
\def\Fin{\mathop{\rm Fin}\nolimits}
\def\br{\mathop{\rm br}\nolimits}
\def\fin{\mathop{\rm fin}\nolimits}
\def\Ann{\mathop{\rm Ann}\nolimits}
\def\mspec{\mathop{\rm mspec}\nolimits}
\def\spec{\mathop{\rm spec}\nolimits}
\def\End{\mathop{\rm End}\nolimits}
\def\Mod{\mathop{\rm Mod}\nolimits}
\def\bfb{\mathop{\rm\bf b}\nolimits}
\def\bfd{\mathop{\rm\bf d}\nolimits}
\def\bfi{\mathop{\rm\bf i}\nolimits}
\def\bfj{\mathop{\rm\bf j}\nolimits}
\def\tr{{\rm tr}}
\def\df{{\rm df}}
\def\bfk{\mathop{\rm\bf k}\nolimits}
\def\bfc{\mathop{\rm\bf c}\nolimits}
\def\bEnd{\mathop{\rm\bf End}\nolimits}
\def\id{\mathop{\rm id}\nolimits}
\def\Ass{\mathop{\rm Ass}\nolimits}
\def\Ext{\mathop{\rm Ext}\nolimits}
\def\Mult{\mathop{\rm Mult}\nolimits}
\def\Ines{\mathop{\rm Ines}\nolimits}
\def\Hom{\mathop{\rm Hom}\nolimits}
\def\Mono{\mathop{\rm Mono}\nolimits}
\def\Aut{\mathop{\rm Aut}\nolimits}
\def\iso{\mathop{\rm Iso}\nolimits}
\def\bHom{\mathop{\rm\bf Hom}\nolimits}
\def\Rk{ R_\k-\mathop{\bf Mod}}
\def\Rn{ R_n-\mathop{\bf Mod}}
\def\map{\mathop{\rm map}\nolimits}
\def\cf{\mathop{\rm cf}\nolimits}
\def\op{\mathop{\rm op}\nolimits}
\def\msp{\mathop{\rm msp}\nolimits}
\def\top{\mathop{\rm top}\nolimits}
\def\product{\mathop{\rm prod}\nolimits}
\def\Ker{\mathop{\rm Ker}\nolimits}
\def\Bext{\mathop{\rm Bext}\nolimits}
\def\Br{\mathop{\rm Br}\nolimits}
\def\dom{\mathop{\rm Dom}\nolimits}
\def\min{\mathop{\rm min}\nolimits}
\def\im{\mathop{\rm Im}\nolimits}
\def\max{\mathop{\rm max}\nolimits}
\def\AX{A[X]}
\def\rk{\mathop{\rm rk}}
\def\Diam{\diamondsuit}
\def\E{{\mathbb E}}
\def\N{{\mathbb N}}
\def\Z{{\mathbb Z}}
\def\Q{{\mathbb Q}}
\def\M{{\mathbb M}}
\def\B{{\mathbb B}}
\def\bQ{{\bf Q}}
\def\bF{{\bf F}}
\def\bX{{\bf X}}
\def\bs\boldsymbol{\sigma}
\def\bS{{\mathbb S}}
\def\AA{{\cal A}}
\def\BB{{\cal B}}
\def\CC{{\cal C}}
\def\DD{{\cal D}}
\def\MM{{\cal M}}
\def\TT{{\cal T}}
\def\FF{{\cal F}}
\def\GG{{\cal G}}
\def\PP{{\cal P}}
\def\SS{{\cal S}}
\def\XX{{\cal X}}
\def\YY{{\cal Y}}
\def\fAb{\mathfrak{ Ab}}
\def\fA{\mathfrak{ A}}
\def\fS{{\mathfrak S}}
\def\fH{{\mathfrak H}}
\def\fU{{\mathfrak U}}
\def\fW{{\mathfrak W}}
\def\fK{{\mathfrak K}}
\def\fp{{\mathfrak p}}
\def\fb{{\mathfrak b}}
\def\fc{{\mathfrak c}}
\def\PT{{\mathfrak{PT}}}
\def\T{{\mathfrak{T}}}
\def\fX{{\mathfrak X}}
\def\fP{{\mathfrak P}}
\def\X{{\mathfrak X}}
\def\Y{{\mathfrak Y}}
\def\F{{\mathfrak F}}
\def\C{{\mathfrak C}}
\def\fB{{\mathfrak B}}
\def\J{{\mathfrak J}}
\def\fN{{\mathfrak N}}
\def\fM{{\mathfrak M}}
\def\fz{{\mathfrak z}}
\def\Fk{{\F_\k}}
\def\bar{\overline }
\def\Bbar{\bar B}
\def\Cbar{\bar C}
\def\Pbar{\bar P}
\def\Tbar{\bar T}
\def\Dbar{\bar D}
\def\abar{\bar a}
\def\bbar{\bar b}
\def\fbar{\bar f}
\def\mbar{\bar m}
\def\vbar{\bar \va}
\def\xbar{\bar x}
\def\dbar{\bar d}
\def\ybar{\bar y}
\def\zbar{\bar z}
\def\nubar{\bar \nu}
\def\Abar{\bar A}
\def\a{\alpha}
\def\e{\varepsilon}
\def\f{\phi}
\def\o{\omega}
\def\k{\kappa}
\def\hv{\widehat\va}
\def\va{\varphi}
\def\z{\zeta}
\def\s{\sigma}
\def\l{\lambda}
\def\m{\mu}
\def\lo{\l^{\aln}}
\def\b{\beta}
\def\g{\gamma}
\def\d{\delta}
\def\ale{\aleph_1}
\def\aln{\aleph_0}
\def\Cont{2^{\aln}}
\def\nld{{}^{ n \downarrow }\l}
\def\n+1d{{}^{ n+1 \downarrow }\l}
\def\hsupp#1{[[\,#1\,]]}
\def\size#1{\left|\,#1\,\right|}
\def\Binfhat{\widehat {B_{\infty}}}
\def\Zhat{\widehat \Z}
\def\Mhat{\widehat M}
\def\shat{\widehat \s}
\def\Rhat{\widehat R}
\def\Phat{\widehat P}
\def\Fhat{\widehat F}
\def\fhat{\widehat f}
\def\Ahat{\widehat A}
\def\Chat{\widehat C}
\def\Ghat{\widehat G}
\def\Bhat{\widehat B}
\def\Btilde{\widetilde B}
\def\Ftilde{\widetilde F}
\def\restr{\mathop{\upharpoonright}}
\def\to{\rightarrow}
\def\arr{\longrightarrow}
\newcommand{\norm}[1]{\text{$\parallel\! #1 \!\parallel$}}
\newcommand{\supp}[1]{\text{$\left[ \, #1\, \right]$}}
\def\set#1{\left\{\,#1\,\right\}}
\newcommand{\mb}{\mathbf}
\newcommand{\cp}{\widehat}
\newcommand{\wt}{\widetilde}
\newcommand{\dsum}{\bigoplus}
\newcommand{\norma}[1]{\mbox{$\parallel\! #1 \!\parallel_A$}}
\newcommand{\suppx}[1]{\mbox{$\left[ \, #1\, \right]_X$}}
\newcommand{\suppa}[1]{\mbox{$\left[ \, #1\, \right]_A$}}
\newcommand{\suppl}[1]{\mbox{$\left[ \, #1\, \right]_\l$}}
\newcommand{\card}[1]{\mbox{$\left| #1 \right|$}}
\newcommand{\union}{\bigcup}
\newcommand{\inters}{\bigcap}
\newcommand{\pure}{\subseteq_\ast}
\newcommand{\mapr}[1]{\xrightarrow{#1}}

\title{{\sc Generalized $E$-Algebras via $\l$-Calculus I}
\footnotetext{Supported by the project No. I-706-54.6/2001 of the
German-Israeli Foundation for\\ Scientific Research \& Development\\
GbSh867 in Shelah's archive\\
subject classification (2000):\\
primary: 20K20, 20K30;  \\
secondary: 16S60, 16W20;\\
Key words and phrases: $E$-rings, endomorphism rings} }

\author{ R\"udiger G\"obel and Saharon Shelah}
\date{}
\maketitle
\begin{abstract}
An $R$-algebra $A$ is called an $E(R)$--algebra if the canonical
homomorphism from $A$ to the endomorphism algebra $\End_RA$ of the
$R$-module ${}_RA$, taking any $a \in A$ to the right multiplication
$a_r\in \End_RA$ by $a$, is an isomorphism of algebras. In this case
${}_RA$ is called an $E(R)$--module. There is a proper class of
examples constructed in \cite{DMV}. $E(R)$-algebras arise naturally
in various topics of algebra. So it is not surprising that they were
investigated thoroughly in the last decade, see
\cite{DG2,DV,Fa,FHR,GG, GSS, GSS1,GSt,Pi,PV}. Despite some efforts
(\cite{GSS1, DV}) it remained an open question whether proper
generalized $E(R)$-algebras exist. These are $R$--algebras $A$
isomorphic to $\End_RA$ but not under the above canonical
isomorphism, so not $E(R)$--algebras. This question was raised about
30 years ago (for $R=\Z$) by Schultz \cite{Sch} (see also
Vinsonhaler \cite{Vi}). It originates from Problem 45 in Fuchs
\cite{Fu0}, that asks one to characterize the rings $A$ for which
$A\cong \End_\Z A$ (as rings). We will answer Schultz's question,
thus contributing a large class of rings for Fuchs' Problem 45 which
are not $E$-rings. Let $R$ be a commutative ring with an element
$p\in R$ such that the additive group $R^+$ is $p$-torsion-free and
$p$-reduced (equivalently $p$ is not a zero-divisor and
$\bigcap_{n\in\o} p^nR=0$). As explained in the introduction we
assume that either $\size{R} <\Cont$ or that $R^+$ is free, see
Definition \ref{sigmaS}.

The main tool is an interesting connection between $\l$-calculus
(used in theoretical computer sciences) and algebra. It seems
reasonable to divide the work into two parts; in this paper we
will work in V=L (G\"odels universe) where stronger combinatorial
methods make the final arguments more transparent. The proof based
entirely on ordinary set theory (the axioms of ZFC) will appear in
a subsequent paper \cite{GS}. However the general strategy will be
the same, but the combinatorial arguments will utilize a
prediction principle that holds under ZFC.
\end{abstract}

\section{Introduction to generalized
$\boldsymbol{ E(R)}$-Algebras}\label{secintro}

Let $\bS$ be a countable, multiplicatively closed subset of a
commutative ring $R$ with $1$.  An $R$--module $M$ is
$\bS$--reduced if $\bigcap_{s\in \bS} sM = 0$ and it is
$\bS$--torsion--free if $sm = 0, \ m\in M, s\in \bS$ implies
$m=0$. Suppose that $R$ (as an $R$-module) is $\bS$--reduced and
$\bS$-torsion-free. Then $R$ is called an $\bS$-ring, see
\cite{GT}. In order to avoid zero--divisors as in the case of
$\Z$--adic completion $\prod_p J_p $ of $\Z$ we also assume that
$\bS$ is cyclically generated, i.e. $\bS = \langle p\rangle := \{
p^n : n\in\o\}$ for some $p\in R$. We will concentrate on
$\bS$--cotorsion--free modules. An $\bS$--torsion-free and
$\bS$--reduced $R$-module $M$ is $\bS$--cotorsion--free if
$\Hom(\Rhat,M)=0$, where $\Rhat$ denotes the $\bS$-completion of
$R$. A submodule $U\subseteq M$ is $\bS$--pure (we also write
$U\subseteq_* M$) if $sM\cap U \subseteq sU$ for all $s\in\bS$.
Note that $R$, being $\bS$--reduced, is Hausdorff in the
$\bS$--topology. In the proof of Step Lemma \ref{step3} we will
also use the following condition on the additive group $R^+$ of
$R$ which implies that $R$ is $\bS$-cotorsion-free.

\begin{definition}\label{sigmaS} An $R$-module $M$ is
{\em $\Sigma\bS$-incomplete} if for any sequence $0\ne m_n\in M$
($n\in \o$) there are $a_n\in \{0,1\}$ with
$\sum_{n\in\o}p^na_nm_n\notin M$. If $M=R^+$ we say that $R$ is
$\Sigma\bS$-incomplete.
\end{definition}

All $\bS$-rings of size $<\Cont$ are $\Sigma\bS$-incomplete as shown
in \cite{GM}. Thus it follows easily that any $\bS$-ring which is a
direct sum of $\bS$-invariant subgroups of size $<\Cont$ is
$\Sigma\bS$-incomplete as well. So we deduce from \cite{GM} a

\begin{corollary}\label{corGM} If an $\bS$-ring $R$ is a direct sum of
$\bS$-invariant subgroups of size $<\Cont$, then $R$ is
$\Sigma\bS$-incomplete. In particular, if $\ \bS$ generates the
ordinary $p$-adic topology (i.e. for $1\in R$ there is $p\in
\langle 1\rangle$ and $\bS = \langle p\rangle$) and the additive
group $R^+$ is free, then $R$ is $\Sigma\bS$-incomplete.
\end{corollary}

We recall the main definition.

\begin{definition} If $A$ is an $R$--algebra, then $\delta: A \arr End_RA$
denotes the homomorphism which takes any $a \in A$ to the
$R$--endomorphism $a\d=a_r$ which is multiplication by $a$ on the
right. If this homomorphism is an isomorphism, then $A$ is called an
$E(R)$--algebra and ${}_RA$ is called an $E(R)$--module. By ${}_RA$
we denote the $R$-module structure of an $R$-algebra $A$.
\end{definition}

$E(R)$-algebras can also be defined dually, assuming that the
homomorphism $$\End_RA\arr A \ (\va \arr 1\va)$$ is an isomorphism.
It is easy to see that $E(R)$--algebras are necessarily commutative.

For any $\bS$-ring (with $\bS$ cyclically generated) that is
$\Sigma\bS$-incomplete we will construct non--commutative
$R$--algebras with $\End_RA\cong A$. Hence these $A$s are
generalized $E(R)$-algebras but not $E(R)$--algebras. If $R=\Z$ and
${}_RA$ is an abelian group, then we do not mention the ring $\Z$:
e.g. $E(\Z)$--modules are just $\E$--groups. The existence of
generalized $\E$--rings answers a problem in \cite{Sch,Vi}.

If $\k$ is a cardinal, then let $\k^o=\{\a : \cf(\a)=\o, \a\in\k\}$.
We will only need the existence of a non-reflecting subset
$E\subseteq \k^o$ for some regular uncountable, not weakly compact
cardinal $\k$ such that the diamond principle $\Diam_\k E$ holds. It
is well-known (see e.g. Eklof, Mekler \cite{EM1}) that $\Diam_\k E$
is a consequence for all non-reflecting subsets $E$ of regular
uncountable, not weakly compact cardinals $\k$ in G\"odel's universe
(V = L). We indicate our (weaker) set theoretic hypothesis (which
also holds in other universes) as $\Diam_\k E$ in our following main
result.

\begin{theorem}\label{mainth} Let $R$ be a $\Sigma\bS$-incomplete
$\bS$-ring for some cyclically generated $\bS$. If $\k>\size{R}$
is a regular, uncountable cardinal and $E\subseteq \k^o$ a
non-reflecting subset with $\Diam_\k E$, then there is an
$\bS$--cotorsion--free, non--commutative $R$--algebra $A$ of
cardinality $|A| = \k $ with $\End_RA \cong A$. Moreover any
subset of cardinality $<\k$ is contained in an $R$-monoid-algebra
of cardinality $<\k$.
\end{theorem}

A similar result without the set theoretic assumption will be
shown in \cite{GS}:

\begin{theorem} Let $R$ be a $\Sigma\bS$-incomplete
$\bS$-ring for some cyclically generated $\bS$. For any cardinal
$\k  = \mu^+$ with $\size{R}\le \mu^{\aleph_0}=\mu$ there is an
$\bS$--cotorsion--free, non--commutative $R$--algebra $A$ of
cardinality $|A| = \k $ with $\End_RA\cong A$.
\end{theorem}

It seems particularly interesting to note that the $R$-monoid $A$
comes from (classical) $\l$-calculus taking into account that
elements of an $E(R)$-algebra $A$ are at the same time
endomorphisms of $A$, thus the same phenomenon appears as known
for computer science and studied intensively in logic in the
thirties of the last century. The problems concerning the
semantics of computer science were solved four decades later by
Scott \cite{Sc1,Sc2}. We will describe the construction of the
underlying monoid $M$ explicitly. Since this paper should be
readable for algebraists with only basic background on model
theory, we will also elaborate the needed details coming from
model theory. The basic knowledge on model theory is in \cite{Ro},
for example. In Section 4 and 5 the monoid $M$ will be completed
and become the algebra $A$.

\section{Model theory of bodies and skeletons via $\boldsymbol{\l}$-calculus}\label{secskel}

\subsection{Discussion}

We begin by defining terms for a skeleton and will establish a
connection with $\l$-calculus. Let $R$ be any commutative $\bS$-ring
with $\bS=\langle p\rangle$. ($\Sigma\bS$-incompleteness will be
added in Section \ref{S5}.)

By definition of generalized $E(R)$-algebras $A$, endomorphisms of
${}_RA$ must be considered as members of $A$. Hence they act on
${}_RA$ as endomorphisms while they are elements of ${}_RA$ at the
same time. Thus we will introduce the classical definitions from
$\l$-calculus over an infinite set $X$ of {\em free} variables and
an infinite set $Y$ of {\em bound} variables to represent those
maps. First note that we can restrict ourselves to unary, linear
functions because endomorphisms are of this kind. (The general
argument to reduce $\l$-calculus to unary functions was observed
by Sch\"onfinkel, see \cite[ p. 6]{Ba}.) What are the typical
terms of our final objects, the bodies? If $x_1$ and $x_2$ are
members of the generalized $E(R)$-algebras $A$ and $a,b\in R$,
then also polynomials like $\s_n(x_1,x_2) = a x_1^n + b x_2^3$
belong to the algebra $A$, so there are legitimate functions
$p_n(y)=\l y.\s_n(y,x_2)$ on $A$ taking $y\arr \s_n(y,x_2)$ and
$A$ must be closed under such \lq generalized polynomials'. This
observation will be described in Definition \ref{reppol} and taken
care of in Proposition \ref{linterm} and in our Main Lemma
\ref{mainl}. A first description of these generalized polynomials
will also be the starting point for our construction and we begin
with its basic settings.

\subsection{The notion of terms}

Let $\tau$ be a vocabulary with no predicates; thus $\tau$ is a
collection of function symbols with an arity function $\tau\arr \o$
defining the places of function symbols. Moreover let $X$ be an
infinite set of free variables. Then unspecified ($\tau,X$)-terms
(briefly called `terms') are defined inductively as the closure of
the atomic terms under these function symbols (only), that is:
\begin{enumerate}
    \item Atomic terms are the $0$-place functions: the
    individual constants (in our case $1$) and members $x$ from
    $X$.
    \item The closure: If $\s_0,\dots,\s_{n-1}$ are terms and
    $F$ is an $n$-place function symbol from $\tau$, then
    $F(\s_0,\dots,\s_{n-1})$ is a term.
\end{enumerate}

We also define the (usual) length $l(\s)$ of a term $\s$
inductively: Let $l(\s)=0$ if $\s$ is atomic and $l(\s)=k+1$ if
$\s$ derives from (ii) with $k=\max\{l(\s_i): i<n\}$.

If $\s$ is an unspecified $(\tau,X)$-term, then we define (also by
induction on $l(\s)$) a finite subset $FV(\s)\subset X$ of {\em
free variables of $\s$}:
\begin{enumerate}
    \item [(a)] If $\s$ is an individual constant, then
    $FV(\s)=\emptyset$ and if $\s\in X$, then $FV(\s)=\{\s\}$.
    \item [(b)] If $\s= F(\s_0,\dots,\s_{n-1})$ is defined as in
    (ii), then $FV(\s)=\bigcup_{i<n}FV(\s_i)$.
\end{enumerate}

We fix some further notation. Let $\xbar=\langle
x_0,\dots,x_{n-1}\rangle$ be a finite sequence of members $x_i$ from
$X$ without repetitions and $\Im(\xbar)=\{x_i : i<n\}$ (and
similarly $\xbar'$). Then we define the (specified)
$(\tau,X)$-terms: A $(\tau,X)$-term is a pair $\boldsymbol{\sigma}=
(\s,\xbar)$ with $\s$ an unspecified $(\tau,X)$-term and $\xbar$ a
finite sequence from $X$ with $FV(\s)\subseteq \Im(\xbar)$. If there
is no danger of confusion, then we will also write
$\boldsymbol{\sigma}= \s(\xbar)$. If $t(\tau,X)$ is the set of all
$(\tau,X)$-terms, then $t(\tau):=\{\s: \boldsymbol{\sigma} \in
t(\tau,X)\}$ is the set of all unspecified $(\tau,X)$-terms.
Furthermore observe that for $\xbar\subseteq \xbar'$ (as maps) with
$(\s,\xbar)$ a $(\tau,X)$--term also $(\s,\xbar')$ is a
$(\tau,X)$--term. For $(\tau,X)$--terms we can define a natural
substitution: if $(\s,\xbar)\in t(\tau,X)$ with $\xbar=\langle
x_0,\dots,x_{n-1}\rangle$ and $\s_0,\dots,\s_{n-1}\in t(\tau)$, then
substitution is defined by
 $$\mbox{Sub}\,^{\langle x_0,\dots,x_{n-1}\rangle}_{\langle
\s_0,\dots,\s_{n-1}\rangle} (\s,\xbar):= \s(\s_0,\dots,\s_{n-1}),$$
 replacing every occurrence of $x_i$ by $\s_i$. If we replace
(if necessary) free variables of the $\s_i$, we can find a sequence
$\xbar'$ with $(\s(\s_0,\dots,\s_{n-1}),\xbar')\in t(\tau,X)$. This
is a good place for two standard notations: let $\bbar=\langle
b_0,\dots, b_{n'-1}\rangle$ be a finite sequence of elements without
repetition from a set $B$. If $n=n'$ and if $(\s,\xbar)$ is as
above, we say that $\bbar$ is a sequence from $B$ ( suitable)  for
$\xbar$ and write $\s(\bbar) =
\mbox{Sub}\,^{\xbar}_{\bbar}(\s,\xbar)$.

A free variable $x\in X$ is a  dummy variable of $(\s,\xbar)$, if
$x\in \Im(\xbar)\setminus FV(\s)$, and we say that $(\s,\xbar)$ is
$X$--reduced, if it has no dummy variables, i.e. $FV(\s)=
\Im(\xbar)$. Trivially, for any $(\s,\xbar)$ we get a natural
$X$--reduced term  by removing those entries of $\xbar$ that
correspond to dummy variables. In this case $\xbar =\langle
x_0,\dots,x_{n-1}\rangle$ becomes $\xbar'=\langle x_{i_0}, \dots,
x_{i_{t-1}}\rangle$ for some $0\le i_0 < i_1 < \dots <i_{t-1}\le
n-1$ and we can use substitution to replace $\xbar'$ by the more
natural sequence $\xbar''=\langle x_0,\dots,x_{t-1}\rangle$: if
$\s':=\mbox{Sub}\,^{\langle x_{i_0}, \dots,
x_{i_{t-1}}\rangle}_{\langle x_0,\dots,x_{t-1}\rangle}(\s,\xbar')$,
then $(\s,\xbar')=(\s',\xbar'')$ (by an axiom below).

\subsection{The vocabulary of a skeleton and its laws}
Let $Y$ be an infinite set of so-called {\em bound} variables (used
as variables for function symbols) and (as before) let
$\ybar=\langle y_0,\dots,y_{n-1}\rangle$ be a finite sequence of
elements from $Y$ without repetitions. Also in this particular case
of the {\em vocabulary $\tau^{sk}$ of a skeleton} the collection
$\tau^{sk}$ will consists of an individual constant $1$, of
variables and of function symbols (only), defined inductively as
$\tau^{sk}_k\ (k\in \o)$; moreover, let $\tau^{sk}_{<k}:=
\bigcup_{m<k}\tau^{sk}_m$ for $k\le\o$, $\tau^{sk}_{\le k}:=
\bigcup_{m\le k}\tau^{sk}_m$ for $k<\o$ and $\tau^{sk}
:=\tau^{sk}_{<\o}$.
\begin{enumerate}
    \item  (Step k=0) The vocabulary $\tau^{sk}_0$ consists of an
    individual constant $1$, free variables $x\in X$ and bound
    variables $y\in Y$. Moreover, we need a particular `$2$--arity word
    product' function symbol $F_\odot$ such that
    $F_\odot(x_0,x_1) = x_0x_1$ is {\it concatenation}.
    \item  (Step $k=m+1$) Suppose that $\tau^{sk}_{\le m}$ is defined
    and $(\s,\xbar)\in t(\tau^{sk}_{\le m},X)$ is a specified
    term of length $k$ with
    $\xbar=\langle x_0,\dots,x_{n-1}\rangle$ and
    $\ybar=\langle y_0,\dots,y_{n-1}\rangle$ suitable for $\xbar$, then
    $F_{\s(y_0,\dots,y_{n-1})}$ is an $n$--place function symbol
    belonging to $\tau^{sk}_k$ (but not to $\tau^{sk}_{<k}$).
\end{enumerate}
 For the collection of terms of the skeletons we will
write $t(\tau^{sk}_{<k},X)$ (where $k\le\o$ is as above) and
$t(\tau^{sk},X):=t(\tau^{sk}_{<\o},X)$. Its members $(\s,\xbar)$
will also be called {\em (generalized) monomials} (because they are
expressed as products).

We now define inductively

\begin{eqnarray} \label{theoskeleton}\text{{\bf The theory of skeletons,
i.e. the axioms $T^{sk}_{<k}$ for $\boldsymbol{\tau^{sk}_{< k}\
(k\le \o)}$:}}\end{eqnarray}

\bigskip
In the following let $\xbar =\langle x_0,\dots,x_n\rangle$, put
$T^{sk}_{<k}:= \bigcup_{m<k}T^{sk}_m$ and $T^{sk}:=T^{sk}_{<\o}$.

\begin{enumerate}
    \item(Step k=0) If $x\in X$, then $1x = x1 = x$ and $1\cdot 1=1$ belong to $T^{sk}_0$.
    \item(Step k=m+1) $T^{sk}_k$ comprises the following laws:
    \begin{enumerate}
        \item If $(\s,\xbar)\in t(\tau^{sk},X)$, $x_0\in FV(\s)$ and $F_{\s(y_0,\dots,y_n)}$
        is a function symbol in $\tau^{sk}_{\le k}\setminus \tau^{sk}_{<k}$
        related to the term $(\s,\xbar)$, then
        $$xF_{\s(y_0,\dots,y_n)}(x_1,\dots,x_n) =
        \s(x,x_1,\dots,x_n).$$
        \item If $(\s,\xbar),(\s,\xbar')$ are $(\tau^{sk}_{\le
        k},X)$-terms with $\xbar\subseteq \xbar':=\langle
        x_0,\dots,x_{n'}\rangle$, then
        $$ F_{\s(y_0,\dots,y_n)}(x_1,\dots,x_n)=F_{\s(y_0,\dots,y_{n'})}(x_1,\dots,x_{n'}).$$
        \item If $\pi$ is an injective map $\{1,\dots,n\}\arr \o\setminus \{0\}$ and
        $\s'(x_0,x_1,\dots,x_n) := \s(x_0,x_{\pi(1)},\dots,x_{\pi(n)})$, then
        $$F_{\s(y_0,\dots,y_n)}(x_1,\dots,x_n)=F_{\s'(y_0,\dots,y_n)}
        (x_{\pi(1)},\dots,x_{\pi(n)}).$$
        \item If $(\s_i,\xbar)\in \tau^{sk}_{\le m}$ for $i=1,2$ and $T^{sk}_{\le
        m}\vdash (\s_1,\xbar)= (\s_2,\xbar)$, then
        $$F_{\s_1(y_0,\dots,y_n)}(x_1,\dots,x_n)=F_{\s_2(y_0,\dots,y_n)}(x_1,\dots,x_n).$$
    \end{enumerate}
\end{enumerate}

\begin{remark} \label{l-def} \begin{enumerate}
 \item Recall that $T \vdash (\s,\xbar)$, means that $(\s,\xbar)$
 follows by the axioms $T$. For convenience (as for free groups) we denote the
 empty product by $1$.
 \item Using the notion of $\l$-calculus for (ii)(a), the unary function
 $F_{\s(y_0,\dots,y_n)}(x_1,\dots,x_n)$ is $\l y_0.\s(y_0,x_1,\dots,x_n)$
 and it acts as $x\l y_0.\s(y_0,x_1,\dots,x_n)= \s(x,x_1,\dots,x_n).$
\end{enumerate}
\end{remark}

The axioms in $T^{sk}_{<k} \ (k\le \o)$ are equations; thus we
have an immediate important application from varieties.

\begin{observation}\label{theories} The theories $T^{sk}_{<k} \ (k\le \o)$  are
varieties with vocabulary $\tau^{sk}_{<k}$. A model $M$ of
$T^{sk}_{<k}$ is an algebra satisfying the axioms of $T^{sk}_{<k}$
and there are models generated freely by any given set.
\end{observation}

\Proof See Gr\"atzer \cite[p. 167]{Gr} or Bergman \cite[Chapter
8]{Be}. \fine

We immediately derive  one of our central definitions.

\begin{definition} Let $T^{sk}:= T^{sk}_{<\o}$ and $\tau^{sk}
=\tau^{sk}_{<\o}$ be as in Observation \ref{theories}. Any
$T^{sk}$-model (an algebra satisfying $T^{sk}$) is called a {\em
skeleton} and two skeletons are called isomorphic if they are
isomorphic as $T^{sk}$-models; see e.g. \cite[p. 262]{Be} or
\cite[p. 5]{Ro}.\end{definition}

For applications it is useful to recall the following
\begin{definition}\label{freebasis} of (free) generators of a $T^{sk}_{<k}$-skeleton.
\begin{enumerate}
    \item The $T^{sk}_{<k}$-model $M$ is {\em generated} by a set
    $B\subseteq M$, if for any $m\in M$ there are $(\s,\xbar)\in
    t(\tau^{sk}_{<k},X)$ and a sequence $\bbar$ in $B$ suitable for $\xbar$ with
    $F_{\s(\ybar)}\bbar= m$.
    \item The $T^{sk}_{<k}$-model $M$ is {\em freely generated} by $B\subseteq
    M$ ($B$ is a basis of $M$) if $B$ generates $M$ and if for any $\bbar$ and $(\s,\xbar),
    (\s',\xbar)\in t(\tau^{sk}_{<k},X)$ with $\bbar$  from $B$ suitable for $\xbar$,
    and $F_{\s(\ybar)}\bbar= F_{\s'(\ybar)}\bbar$, follows
    $(\s,\xbar) =(\s',\xbar)$ from $T^{sk}_{<k}$.
\end{enumerate}
\end{definition}

For any set $B$ we will construct a skeleton $\B$ freely generated
by $B$. For this we need

\subsection{Reduction of terms}
 Freeness can easily be checked by the
usual rewriting process (as in group theory). Thus we define for
each term $(\s,\xbar)\in t(\tau^{sk}_{<k},X)$ its reduced form
$red(\s,\xbar):=(\s^r,\xbar^r)\in t(\tau^{sk}_{<k},X)$.
Inductively we apply the axioms (\ref{theoskeleton}) (in
particular (ii)(a) which connect formulas with function symbols)
to shorten the length of a term; note that by the axioms
(\ref{theoskeleton}) terms remain the same; we arrive at an
essentially unique reduced term. We first consider the reduction
of unspecified terms and find $\s^r$ from $\s$:

\begin{eqnarray}\label{reduce} { \text{\bf The reduction of terms:}}\end{eqnarray}
\begin{enumerate}
    \item If $\s =x \in X$, then $\s^r= x$ and if $\s=1$, then
    $\s^r=1$.
    \item If $\s=\s'\s''$ and $\s', \s''$ are reduced, but $\s''$ is
    not of the form $F_{\s_0(y_0,\dots,y_{t})}(\s_1,\dots,\s_t)$, then
    $\s^r= \s$ is reduced.
    \item Suppose that $\s'$ and $\s_i \ (i\le t)$ are reduced and
    $\s=F_{\s'(y_0,\dots,y_{t})}(\s_1,\dots,\s_t)$,
    with $(\s',\langle x_0,\dots x_t\rangle)\in t(\tau^{sk}_{<k},X)$ the corresponding
    specified term. First we get rid of dummy variables: let $u:=\{ i\in
    \{1,\dots,t\}, x_i\in FV(\s')\}$; say $u=\{1\le i_1< \dots
    <i_{|u|}\le n\}$ and $\xbar'=\langle x_i : i\in u\rangle$.
    Then $\s^r= F_{\s'(y_0,y_{i_1},\dots,
    y_{i_{|u|}})} \langle \s_i: i\in u\rangle$. See below for a normalization.
    \item If $\s=\s'\s''$ and $\s',\s''$ are reduced terms, but $\s''$ is
    a unary function of the form $F_{\s_0(y_0,\dots,y_{t})}(\s_1,\dots,\s_t)$ then
    $\s^r=\s_0(\s',\s_1,\dots,\s_t)$.
\end{enumerate}

We are ready for the

\begin{definition}\label{reducedterm} An unspecified term $\s\in
t(\tau^{sk}_{<k})$ is {\em reduced} if $\s^r =\s$. A term
$(\s,\xbar)\in t(\tau^{sk}_{<k},X)$ is {\em reduced} if $\s^r =\s$
and $(\s,\xbar)$ has no dummy variables, i.e. $FV(\s) =\Im(\xbar)$.
Thus $red(\s,\xbar)=(\s,\xbar)$. Moreover, $(\s,\xbar)$ is {\em
normalized} if the free variables of $\s$ are enumerated as $\langle
x_0,\dots,x_{n-1}\rangle$. \end{definition}

It is now easy to extend the reduction inductively to
$t(\tau^{sk}_{<k},X)$. Let $red(x,\xbar) =(x,\langle x\rangle)$ and
$red(1,\xbar) = (1,\emptyset)$. In (ii) we first ensure (by free
substitution) that the free variables of $\s'$ and $\s''$ are
disjoint and then order their union. We want to normalize
$(\s^r,\xbar')$ in (iii): the sequence $\xbar'$ is of the form
$\langle x_{i_1},\dots, x_{i_{|u|}}\rangle$ for some $1\le i_1 <
\dots < i_{|u|}\le n$. We will replace it by the more natural
sequence $\xbar'' =\langle x_1,\dots, x_{|u|}\rangle$ and use
substitution $\s''=Sub^{\langle x_{i_1},\dots, x_{i_{|u|}}\rangle}
_{\langle{x_1},\dots, x_{|u|}\rangle} (\s^r,\xbar')$, thus
$T^{sk}_{<k}$ implies  $(\s^r,\xbar') = (\s'',\xbar'')$ by
(\ref{theoskeleton})(ii)(b). In (iv) we order the union of the free
variables $FV(\s_i)\ (i\le t)$ after making them pairwise disjoint
by free substitutions.
\bigskip

Thus we have a definition and a consequence of the last
considerations.

\begin{defobs}\label{redequal} Every term $(\s,\xbar)$ can
be reduced to a (normalized) reduced term $red(\s,\xbar)$ with
$T^{sk}_{<k}\vdash red(\s,\xbar)= (\s,\xbar)$. Let
$t^r(\tau^{sk}_{<k},X)$ be the family of reduced terms from
$t(\tau^{sk}_{<k},X)$; moreover let $t^r(\tau^{sk}_{<k})=\{\s:
\boldsymbol{\sigma}\in t^r(\tau^{sk}_{<k},X)\}$.
\end{defobs}

Thus we consider only elements from $t^r(\tau^{sk}_{<k},X)$ (so in
particular function symbols $F_\s$ have attached reduced terms
$(\s,\xbar)$). We want to discuss how much reduced terms can
differ if they represent the same element of a free skeleton. We
first give the definition which describes this.

\begin{definition} \label{defessequal} Using induction, we say when two reduced elements
$\s_1,\s_2\in t^r(\tau^{sk}_{<k})$ are {\em essentially equal} and
will write $\s_1\doteq \s_2$.
\begin{enumerate}
    \item If $\s_1$ is atomic, then $\s_2$ is atomic and $\s_1=\s_2$.
    \item If $\s_1=\s_1'\s_1''$ such that $\s_1', \s_1''$ are reduced, then $\s_2=\s_2'\s_2''$
    and $\s_1'\doteq\s_2',\s_1''\doteq\s_2''$.
    \item If $\s_i =
    F_{\s'_i(y_0,\dots, y_{m_i})}(\s^i_1,\dots,\s_{m_i}^i)$ for
    $i\le 2$, then $m_1=m_2$ and there is a permutation $\pi$ of
    $\{1,\dots m_1 \}$ with $\s^1_j\doteq\s^2_{\pi(j)}$ for all $j\le m_1$ and also
    $\s'_1\doteq \s'_2$.
\end{enumerate}
\end{definition}

\begin{observation} \label{essequal}
\begin{enumerate}
    \item The relation $\doteq$ is an equivalence relation on
    $t^r(\tau^{sk}_{<k},X)$.
    \item If $F$ is an $n$-place function symbol in $\tau^{sk}_{<k}$ and
    $\s_i\doteq \s_i'\in t^r(\tau^{sk}_{<k},X)$ for $i<n$, then
    $F(\s_0,\dots,\s_{n-1})^r\doteq
    F(\s'_0,\dots,\s'_{n-1})^r$.
\end{enumerate}
\end{observation}

\Proof This is immediate by induction using (\ref{theoskeleton}).
\fine

Using normalization of $\xbar$ and $\bbar$ from
Definition-Observation \ref{redequal} we can deduce a

\begin{proposition} If $M$ is a $T^{sk}_{<k}$-model and $(\s,\xbar)
\in t(\tau^{sk}_{<k},X)$ with $red(\s,\xbar)= (\s^r,\xbar^r)$, then
 $M\vdash \s(\bbar)= \s^r(\bbar^r)$ for any
sequence $\bbar$ from $M$ suitable for $\xbar$ with a sequence
$(\bbar^r)$ obtained by normalization.
\end{proposition}

Note that reduction of terms is defined for each $k\le\o$, thus
formally it depends on $k$. Moreover, $(\tau^{sk}_{<h},X)
\subseteq(\tau^{sk}_{<k},X)$ for all $h\le k\le\o$. Next we show
that the reduction of terms in $(\tau^{sk}_{<h},X)$ is the same even
if it takes place in $(\tau^{sk}_{<k},X)$, i.e.
$t^r(\tau^{sk}_{<h},X)= t^r(\tau^{sk}_{<k},X)\cap
t(\tau^{sk}_{<h},X)$. A similar argument holds for freeness.

\begin{proposition}\label{startred} Let $h\le k\le \o$ and $(\s,\xbar)\in
 t(\tau^{sk}_{<h},X)$.
 \begin{enumerate}
    \item $red(\s,\xbar)$ w.r.t. $t(\tau^{sk}_{<h},X)$ is the same as
    w.r.t. $ t(\tau^{sk}_{<k},X)$.
    \item If $M_z$ is freely generated by $B$ w.r.t.
    $t(\tau^{sk}_{<z},X)$ for $z\in \{h, k\}$, then there is an embedding
    $\iota: M_h\arr M_k$ with $\iota\restr B=\id_B$.
    \item If $\s\in t(\tau^{sk}_{<k})$, then $T^{sk}_{<k}\vdash \s= \s^r$.
 \end{enumerate}
\end{proposition}

\Proof  (i) follows because reduction of elements from
$t(\tau^{sk}_{<h},X)$ only uses terms from $t(\tau^{sk}_{<h},X)$. In
(ii) we can extend the identity $\id:B\arr B$ naturally by induction
to $M_h\arr M_k$ and (iii) follows from Definition-Observation
\ref{redequal}. \fine

The skeleton has the following important property.

\begin{corollary} For $\s_1, \s_2 \in \tau^{sk}$ the following are equivalent
\begin{enumerate}
    \item $T^{sk} \vdash \s_1 = \s_2$.
    \item $\s_1^r\doteq \s^r_2$.
\end{enumerate}
\end{corollary}

\Proof $(ii)\arr (i)$: From the Definition-Observation
\ref{redequal} follows $T^{sk}\vdash \s_1=\s_1^r, \s_2=\s_2^r$ and
by Observation \ref{essequal} is $T^{sk}\vdash \s^r_1=\s_2^r$ thus
$T^{sk} \vdash \s_1 = \s_2$.

$(i)\arr (ii)$: From (i) follows $T^{sk}\vdash \s^r_1=\s_2^r$. Thus
$\s_1^r\doteq \s^r_2$ by Definition \ref{defessequal} and
Observation \ref{essequal}. \fine

\subsection{The skeleton freely generated by X}
Next we construct and discuss free skeletons based on reduced
terms. We will use the infinite set $X$ of free variables to
construct a skeleton $M_X$ which is freely generated by a set
which corresponds by a canonical bijection to  $X$.

By Observation \ref{essequal}(i) we have an equivalence relation
$\doteq$ on $\tau^{sk}: = \tau^{sk}_{<\o}$ with equivalence classes
$[\s]$ for any $\s \in \tau^{sk}$. Let
$$M_X = \{[\s] : \s \in \t(\tau^{sk})\}.$$
The equivalence classes $[\s]$ of atoms $\s$ are singletons by
Definition \ref{defessequal}(i). If $B = \{[x] : x \in X \}$, then
$\iota:X \arr B \ (x \mapsto [x])$ is a bijection and we see that
$M_X$ is a skeleton with basis $B \subseteq M_X$. Moreover, $[\ ]$
is compatible with the application of function symbols:

If $F=F_\s \in \tau^{sk}$ with $(\s,\xbar)\in t(\tau^{sk},X)$ and
$\xbar=\langle x_0,\dots,x_{n-1}\rangle$ is an $n$-place function
symbol and $\s_i\doteq \s_i' \ (i<n)$, then
$F(\s_0,\dots,\s_{n-1})^r\doteq F(\s'_0,\dots,\s'_{n-1})^r$ by
Observation \ref{essequal}(ii), thus
$$F([\s_0],\dots, [\s_{n-1}]) = [(F(\s_0,\dots,\s_{n-1})]\ \
(\s_i \in t^r(\tau^{sk}))$$ is well-defined as follows from
Observation \ref{essequal}(ii).

We have the following

\begin{theorem}\label{thmskeleton} If $X$ is an infinite set (of
free variables) and $M_X$ is defined as above, then the following
holds.\begin{enumerate}
    \item [(a)] $M_X = \{[\s] : \s \in \t^r(\tau^{sk})\}$.
    \item [(b)] $M_X$ is a skeleton with $n$-place functions
    $$[F] : (M_X)^n \arr M_X \ ([\s_0],\dots, [\s_{n-1}]) \mapsto [F(\s_0,\dots,
    \s_{n-1})]$$
    for each $n$-place function symbol $F =
    F_{\s(y_0,\dots,y_{n-1})}$ for $(\s,\xbar) \in
    t^r(\tau^{sk},X)$ with $FV(\s) = \{x_0,\dots,x_{n-1}\}$.
    \item [(c)] $M_X$ is freely generated by $B=\{[x],x\in X\}$, called the
    {\em free skeleton over $X$}. Using $\iota$ above we
    identify $B$ and $X$.
\end{enumerate}
\end{theorem}

\Proof The axioms (\ref{theoskeleton}) are satisfied, e.g. the
crucial condition (ii)(a) follows by definition of $[F]$. \fine

\begin{remark} In the construction of the free skeleton $M_X$ we
also used an infinite set $Y$ of bound variables. However, it
follows by induction that another infinite set $Y'$ of bound
variables leads to an isomorphic copy of $M_X$. Thus we do not
mention $Y$ in Theorem \ref{thmskeleton}. \end{remark}

\begin{lemma} Let $B$ be a subset of the $T^{sk}_{<k}$-model $M$ for $k\le\o$. Then
$B$ is a basis if and only if the following two conditions hold.
\begin{enumerate}
 \item If $c\in M$, then there are $(\s,\xbar)\in
 t^r(\tau^{sk}_{<k},X)$ and a sequence $\bbar$ for $\xbar$ from $B$
 such that $\s(\bbar) =c$.
 \item If $(\s,\xbar), (\s',\xbar')\in t^r(\tau^{sk}_{<k},X)$ with
 $\xbar=\langle x_0,\dots, x_n\rangle, \xbar'=\langle x_0,\dots,
 x_{n'}\rangle$ and $\bbar, \bbar'$ are suitable sequences for
 $\xbar$, $\xbar'$, respectively from $B$, then
 $\s(\bbar)=\s'(\bbar')$ implies $n=n'$ and there is a permutation
 $\pi$ of $\{0,\dots,n\}$ such that $\s(\xbar) =\s'(x_{\pi(0)},
 \dots,x_{\pi(n)})$ and $b'_i= b_{\pi(i)}$ for all $i\le n$.
\end{enumerate}
\end{lemma}

\Proof If $B$ is a basis of $M$, then by Definition
\ref{freebasis} the two conditions of the lemma hold; see
Proposition \ref{startred} (i) and (iii) for (i). Conversely
suppose that (i) and (ii) hold. It is easy to extend inductively a
bijection $B\arr X$ to an an isomorphism between $M$ and the free
skeleton $M_X$ as in Proposition \ref{startred}. Thus $B$ is a
basis. \fine

\subsection{The vocabulary of bodies and their laws}

Recall that $R$ is an $\bS$-ring of size $<\k$ with $\bS=\langle
p\rangle \subseteq R$ as explained in the introduction. Also
recall that $\tau^{sk}_{<k}$ is the vocabulary of skeletons from
the last section, so in particular $\tau^{sk}=\tau_{<\o}^{sk}$
with similar notations for the axioms $T^{sk}_{<k}$.

We now extend the vocabulary $\tau^{sk}_{<k}$ of skeletons to the
vocabulary $\tau^{bd}_{<k}$ of bodies:

Let $\tau^{bd}_{<k}$ comprise all function symbols from
$\tau^{sk}_{<k}$ (so $\tau^{sk}_{<k}\subseteq \tau^{bd}_{<k}$) and
choose additional function symbols:
\begin{enumerate}
    \item [(0)] An individual constant $0$ (for $0$ of an $R$-module)
    \item [(1)] Let $F_+$ be a binary function symbol (in charge
    of addition in $R$-modules). Thus we will write
    $F_+(y_0,y_1) = y_0+y_1$, as usual.
    \item [(2)] For each $a\in R$ let $F_a$ be a unary function
    symbol (for scalar multiplication by $a$ on the left).
    Thus we will write $F_a(y) = ay$, also as usual.
\end{enumerate}

Repeated application of (1) and (2) leads to finite sums like
$\sum_{i=1}^na_iy_i$ and we will write
$\tau^{bd}:=\tau_{<\o}^{bd}$ and call this the {\em vocabulary of
the bodies}. Again terms can be written as $(\s,\xbar)$ with
$\s\in \tau^{bd}$ (or in $\tau^{bd}_{<k}$ as  for skeletons) with
$FV(\s)= \Im(\xbar)$ for reduced terms. The collection of terms of
the bodies will be $t(\tau^{bd}_{<k},X)$ (where $k\le\o$ is as
above). Its members $(\s,\xbar)$ will also be called {\em
(generalized) polynomials}, because we will show (Lemma
\ref{linearexpr}) that they can be expressed as linear
combinations of generalized monomials (terms from
$t(\tau^{sk}_{<k},X)$ or from $t^r(\tau^{sk}_{<k},X)$,
respectively).

As in case of skeletons we now derive the axioms of the bodies in
order to see that they build a variety as well.
\begin{eqnarray}\label{bodytheory}\text {\bf The theory $\boldsymbol{T^{bd}_k}$
of bodies for $\boldsymbol{k\le\o}$:}\end{eqnarray}
\begin{enumerate}
    \item $T^{sk}_k\subseteq T^{bd}_k$.
    \item {\rm Linearity:} If $F\in \tau_{<k}^{sk}$ is an $n$-place
    function symbol, $a_i\in R$ ($i\le t$) and $1\le l\le n$, then
    $$F(x_1,\dots,x_{l-1},\sum_{i=1}^t  a_ix_{li}, x_{l+1},\dots,
    x_n) = \sum_{i=1}^t a_i F(x_1,\dots,x_{l-1}, x_{li}, x_{l+1},\dots,
    x_n).$$
    \item  {\rm The usual module laws:}
    Let $a,b\in R$ and $w,y,z\in
    M$ ($M$ a $T^{bd}_k$-model).
    \begin{enumerate}
        \item $0+y = y$, $z+y=y+z$, $w+(y+z)=(w+y)+z$.
        \item $1y=y$, $a(by) = (ab)y, a(z+y) = az +ay, (a+b)y = ay + by, y + (-1)y =0$.
    \end{enumerate}
\end{enumerate}

\begin{observation} \label{bdtheories} The theories $T_{<k}^{bd} \ (k\le \o)$  are
varieties with vocabulary $\tau_{<k}^{bd}$. A model $M$ of
$T_{<k}^{bd}$ is an algebra satisfying the axioms of $T_{<k}^{bd}$
and there are models generated freely by any given set.
\end{observation}

\Proof See Gr\"atzer \cite[p. 198, Theorem 3]{Gr} or  Bergman
\cite[Chapter 8]{Be}. \fine

\begin{definition}\label{bodyvar} Let $T^{bd}:= T^{bd}_{<\o}$ and $\tau^{bd}:=\tau^{bd}_{<\o}$
be as in Observation \ref{bdtheories}. Any $T^{bd}$-model (an
algebra satisfying $T^{bd}$) is called a {\em body} and two bodies
are isomorphic if they are isomorphic as $T^{bd}$-models; see e.g.
\cite[p.262]{Be} or \cite{Ro}.
\end{definition}

\begin{observation} Any (generalized) $E(A)$-algebra is a body.
\end{observation}

\Pf Generalized $E(A)$-algebras satisfy $\End_RA=A$. Thus any
function symbol $F_\s$ of $\tau^{bd}$ can be interpreted on $A$ as
a function and the axioms (\ref{theoskeleton}) and
(\ref{bodytheory}) hold. \fine

\bigskip

But note that only free bodies arrive from skeletons, see Section
\ref{bodysec}.

\subsection{ Linearity of unary body functions from
$\boldsymbol{ t(\tau_{<k}^{bd},X)}$}

We will first show that terms in $t(\tau_{<k}^{bd},X)$ are linear
combinations of terms in $t(\tau_{<k}^{sk},X)$, thus every
polynomial (in $t(\tau_{<k}^{bd},X)$) is a linear combination of
monomials (in $t(\tau_{<k}^{sk},X)$).

We will show the following

\begin{lemma}\label{linearexpr} Let $\xbar=\langle x_0,\dots x_{m-1}\rangle$ and
$(\s,\xbar)\in t(\tau_{<k}^{bd},X)$. Then there is $\sum_{l<t}
a_l\s_l(\xbar)$ with
\begin{enumerate}
    \item $(\s_l,\xbar)\in t(\tau_{<k}^{sk},X)$ for $l<t$,
    \item $a_l\in R$ for $l<t$,
    \item and $T^{bd}_{<k}\vdash \s =\sum_{l<t}a_l \s_l$.
\end{enumerate}
\end{lemma}

\Proof (We will now suppress the index `$<k$'.)  We let
$(\s,\xbar)\in t(\tau^{bd},X)$ and prove the lemma by induction on
the length of $\s$. If $\s$ is atomic, then $(\s,\xbar)$ is a
monomial and there is nothing to show.

If $F$ is an $m$-place function symbol from $\tau^{bd}$ and
$\s=F(\s_0,\dots,\s_{m-1})$ with $FV(\s_l)\subseteq FV(\s)$, then by
induction hypothesis for $\s_l$ there are polynomials
$\s_l=\sum_{i<t_l}a_{li}\s_{li}(\xbar)$ with $a_{li}\in R$ and
$\s_{li}$ monomials (terms in $t(\tau^{sk},X)$). We substitute these
sums into $F$ and apply axioms (\ref{bodytheory})(ii) (the
linearity) for functions in the theory of bodies. Thus also $\s$ is
as required.

If $\s = F_+(\s_1,\s_2) = \s_1 + \s_2$ comes from (1) and if $\s=
a\s_1$ arrives from (2) the linearity follows from
(\ref{bodytheory})(iii). Thus the lemma is shown. \fine

\begin{lemma}\label{skellin} If $\xbar=\langle x_0,\dots x_m\rangle$
and $(\s,\xbar)$ is a monomial (a term in $t(\tau^{sk},X)$), then
$\s(\sum_{l<t} x_{0l}a_l,x_1,\dots,x_m) =\sum_{l<t}a_l\s(
x_{0l},x_1,\dots,x_m)$
\end{lemma}
\Proof  This is an easy induction on the length of $\s$:

Case 1: If $\s =1$ and $\s=x_0$, then the claim holds trivially.

Case 2: If $\s=F(\s_0,\dots,\s_m)$, then the claim follows from the
axioms (\ref{bodytheory})(ii) of $T^{bd}$. Similarly, if $\s=
F_+(\s_1,\s_2)$ and $\s=F_a(\s_1)$, then the linearity follows by
definition of these functions and induction hypothesis. \fine

Recall the notion from $\l$-calculus in Remark \ref{l-def}(ii).

\begin{proposition} {\rm\bf The weak completeness of bodies.}
\label{linterm} Let $M$ be a body and $(\s,\xbar)$ a polynomial (a
term in $t(\tau^{bd},X)$) with $\xbar=\langle x_0,\dots,x_m\rangle$
and $d_1,\dots,d_m\in M$. Then there is $(\s',\langle
x_1,\dots,x_n\rangle)\in t(\tau^{bd},X)$ and the following holds.
\begin{enumerate}
    \item $M$ is an $R$-module.
    \item The unary function $\l.z\s(z,d_1,\dots,d_m):
 M\arr M \ (z\mapsto \s(z,d_1,\dots,d_m))$ is the
$R$-endomorphism $\l y.y\s'(z,d_1,\dots,d_m)\in \End_R(M)\ (y\arr
y\s'(d_1,\dots,d_m))$.
\end{enumerate}
\end{proposition}

{\bf Remark.} We will show here that there is a function symbol
$(\s',\xbar)\in t(\tau^{bd},X)$ such that $\s(d,d_1,\dots,d_m) =
d\s'(d_1,\dots,d_m)$ for all $d\in M$, see axioms
(\ref{theoskeleton})(ii)(a) of the skeletons.

\Proof By the axioms (\ref{bodytheory}) of $T^{bd}$ it is clear
that $M$ is an $R$-module. It remains to show (ii). Let
$(\s,\xbar')\in t(\tau^{bd},X)$ with $\xbar'=\langle
x_0,\dots,x_m\rangle$ and $\xbar=\langle x_1,\dots,x_m\rangle$. By
Lemma \ref{linearexpr} there are monomials $(\s_l,\xbar')\in
t(\tau^{sk},X)$ and $a_l\in R$ such that
 $$\s =\sum_{l<t} a_l\s_l.$$
 From the construction of $\tau^{sk}$ we also have
function symbols $F_{\s_l(y_0,\dots,y_m)}$ satisfying the axioms
(\ref{theoskeleton}) and (\ref{bodytheory}). Thus
 $$ \s_l(\xbar')=x_0 F_{\s_l(y_0,\dots,y_m)}(\xbar)$$
 and we put
 $$\s'(\xbar)= \sum_{l<t} a_lF_{\s_l(y_0,\dots,y_m)}(\xbar).$$

For (ii) it remains to show $\s(d,d_1,\dots,d_m)
=d\s'(d_1,\dots,d_m)$ for all $d\in M$ which will follow from
$\s(x_0,\xbar) = x_0\s'(\xbar)$. We use the three displayed
formulas and calculate
 $$\s(x_0,\xbar) = \sum_{l<t} a_l\s_l(x_0,\xbar)=\sum_{l<t} a_l(x_0
F_{\s_l(y_0,\dots,y_m)}(\xbar))$$$$= x_0(\sum_{l<t}
a_lF_{\s_l(y_0,\dots,y_m)}(\xbar))=x_0\s'(\xbar).$$ Hence (ii)
follows. \fine

\section{From the skeleton to the body}\label{secbody}

\subsection{The monoid structure of skeletons}\label{product}

Recall from Theorem \ref{thmskeleton} that the skeleton on an
infinite set $X$ of free variables is the set $M_X = \{[\s] : \s \in
\t^r(\tau^{sk})\}$ with $n$-place functions
 $$[F] : (M_X)^n \arr M_X \ ([\s_0],\dots, [\s_{n-1}]) \mapsto [F(\s_0,\dots,
\s_{n-1})]$$
 for each $n$-place function symbol $F = F_{\s(y_0,\dots,y_{n-1})}$
with $(\s,\xbar) \in t^r(\tau^{sk})$ and $FV(\s) =
\{x_0,\dots,x_{n-1}\}$. For simplicity we will also write Roman
letters for the members of $M_X$, e.g. $m=[\s]\in M_X$. The set
$M_X$ has a distinguished element $1=[1]$ and $m1=1m=m$ holds for
all $m\in M_X$ (thus $M_X$ is an applicative structure with $1$). In
order to turn $M_X$ into a monoid, we first represent $M_X$ as a
submonoid of $\Mono(M_X)$, the injective maps on $M_X$, say $\iota:
M_X\arr \Mono(M_X)$:

Let $a=[\s]\in M_X$ and $\s'\in t(\tau^{bd})$. We will use
induction: If $a=[1]$ then $[\s']a =\s'$, if $a=[x]$, then $[\s']a
=\s'x$ and if $a=[F_{\s(y_0,\dots,y_{n-1})}\langle
x_1,\dots,x_{n-1}\rangle]$ is a unary function as above, then
$[\s']a= \s(\s',x_1,\dots,x_{n-1})$ (so $a\iota = [\l y. y\s]$).
Thus $a\iota$ maps any $m=[\s']\in M_X$ to $m(a\iota)= m (\l y.
y\s)= [m\s]\in M_X$ which can be represented by a reduced element
using (\ref{reduce}). If $a\ne b\in M_X$, then $ 1(a\iota)= a\ne b
=1(b\iota)$ thus $\iota:M_X\arr \Mono(M_X)\subseteq M_X$ is an
embedding. We define multiplication of elements $a,b\in M_X$ as
composition of functions $(a\iota)(b\iota) =(a b)\iota$. This is to
say that from $a=[\s], b=[\s']$ we get the product as the
equivalence class of $\l y.((y\s')\s)$. We will write $a\cdot b= ab$
and will often suppress the map $\iota$. From $\Mono(M_X)$ follows
that also $M_X$ is a monoid. Also note that $[x][x']\ne [x'][x]$ for
any free variables $x,x'\in X$. We get an

\begin{observation}\label{freemono} The free skeleton $(M_X,\cdot,1)$ with composition of
functions as product is a non-commutative (associative) monoid with
multiplication defined as above by the action on $M_X$: If
$[\s],[\s']\in M_X$, then $[\s']\cdot[\s]=[\l y.((y\s')\s)]$.
\end{observation}

\subsection{Free Bodies from skeletons}\label{bodysec}

Finally we will associate with any skeleton $M$ its (canonical) body
$\B_RM$: Let $\B_RM$ be the $R$-monoid algebra $RM$ of the monoid
$M$. Moreover any $n$-place function $F:M^n\arr M$ extends uniquely
by linearity to $F:\B_RM^n \arr \B_RM$. We deduce a

\begin{lemma} \label{freeR} If $R$ is a commutative ring as above and
$M$ a skeleton, then the $R$-monoid algebra $\B_RM$ of the monoid
$M$ is a body. If the skeleton $M_X$ is freely generated by $X$,
then also $\B_RM$ is freely generated by $X$ as a body. Moreover
${}_R\B_RM_X= \bigoplus_{m\in M} mR$.
\end{lemma}

\Proof It is easy to see that $\B_RM$ (with the linear $n$-place
functions) is a body. We first claim that $X$, viewed as $\{[x]:
x\in X\}\subseteq \B_RM_X$ is a basis. First apply Lemma
\ref{linearexpr} to the $R$-monoid $\B_RM_X$: Any $(\s,\xbar)\in
t(\tau^{bd},X)$ can be written as a polynomial $\s=\sum_l\s_la_l$
with monomial $(\s_l,\xbar)\in t(\tau^{sk},X)$. Moreover, any
$\s_l$ is viewed as an element of $\Mono(M_X)$, so axiom
(\ref{theoskeleton})(ii)(a) applies and $\s_l$ becomes a product
of elements from $X$. Thus $X$ generates $\B_RM_X$. The monomials
of the skeleton $M$ extend uniquely by linearity to polynomials of
the free $R$-module ${}_R\B_RM_X= \bigoplus_{m\in M} mR$ from its
basis $M$. \fine

We will also need the notion of an extension of bodies.

\begin{definition} Let $\B$ and $\B'$ be two bodies,
then $\B\le\B'$ ($\B'$ extends $\B$) if and only if $\B\subseteq
\B'$ as $R$-algebras and if $(\s,\xbar) \in t(\tau^{bd},X)$ and
$F_\s$ is a function symbol with corresponding unary, $R$-linear
function $F$ of $\B'$, then its natural restriction to $\B$ is the
function for $\B$ corresponding to $F_\s$.
\end{definition}

\begin{example} Let $X\subseteq X'$ be sets of free variables and $\B,\B'$
be the free bodies generated by the free skeletons obtained from $X$
and $X'$, respectively. Then $\B\le \B'$. In this case we say that
$\B'$ is free over $\B$.
\end{example}

\section{The technical tools for the main construction}\label{mainc}

 The endomorphism ring $\End_R\B_RM_X$ of the $R$-module
${}_R\B_RM_X$ has natural elements as endomorphisms, the linear
maps, our (generalized) polynomials interpreted by the terms in
$\s\in t(\tau^{bd},X)$ acting by scalar multiplication on $\B_RM_X$
as shown in Proposition \ref{linterm}(ii). The closure under these
polynomials is dictated by the properties of $E(R)$-algebras. Thus
we would spoil our aim to construct generalized $E(R)$-algebras if
we `lose these $R$-linear maps' on the way.

\begin{definition}\label{defpoly} Let $(\s,\xbar) \in t(\tau^{bd},X)$ with
$\xbar=\langle x_0,\dots,x_n\rangle$. If $\B$ is a body,
$\dbar=\langle d_1,\dots,d_n\rangle$ with $d_1,\dots,d_n\in \B$,
then we call $s_{\dbar}(y)= \l y.\s(y,\dbar)$ the (generalized)
polynomial over $\B$ with coefficients $\dbar$.\end{definition}

Note that $\s_{\dbar}(y)$ is a sum of products of elements $d_i$ and
$y$. Here we must achieve (full) completeness of the final body,
thus showing that any endomorphism is represented.  By a prediction
principle we kill all endomorphisms that are not represented by
$t(\tau^{bd},X)$ - thus the resulting structure will be complete:
Any $R$-endomorphism of an extended body $\B_RM_X$ will be
represented by a polynomial $q(x)$ over $\B_RM_X$, so $\B_RM_X$ is
complete or equivalently an $E(R)$-algebra.

The fact that $\B_RM_X$ is not just the $R$-linear closure (or
$A$-linear closure for some algebra $A$), makes this final task, to
get rid of undesired endomorphisms harder than in case of realizing
algebras as endomorphism algebras (where the closure is not that
floppy).
\begin{definition}\label{reppol} Let $\B$ be a body and $G={}_R\B$.
Then $\va\in \End_RG$ is called \em{ represented } (by $q(y)$) if
there is a generalized polynomial $q(y)$ with coefficients in $\B$
such that $g\va = q(g)$ for all $g\in G$.\end{definition}

If all elements from $\End_RG$ are represented, then $\B$ is a
generalized $E(R)$-algebra.

As for other algebraic structures we have the

\begin{lemma}\label{charbasis} Let $R$ be an $\bS$-ring as above and
$\B$ be a body generated by $B$, then $B$ is a basis of $\B$ if one
of the following equivalent conditions holds.

\begin{enumerate}
 \item If $B'=\B_RM_{X}$ is the body generated by the free skeleton
 $M_{X}$ and $X\arr B$ is a bijection, then this map extends to an
 isomorphism $\B'\arr \B$
 \item  $B$ is independent in $\B$, i.e. if $(\s_1,\xbar),(\s_2,\xbar)\in t^r(\tau^{bd},X)$
 and the sequence $\bbar$ from $B$ is suitable for $\xbar$ such that
 $\s_1(\bbar)=\s_2(\bbar)$, then $T^{bd}\vdash \s_1(\xbar)=\s_2(\xbar)$.
 \item For all bodies $H$ and maps $\va: B\arr H$ there is an extension
 $\vbar: \B_RM_X\arr H$ as $T^{bd}$-homomorphism.
\end{enumerate}
\end{lemma}

\Pf The proof is well known for varieties (see Gr\"atzer \cite[p.
198, Theorem 3]{Gr} or Bergman \cite[Chapter 8]{Be}), so it
follows from Observation \ref{bdtheories}.\fine

\begin{fprop}\label{fproper} Let $R$ be an $\bS$-ring as above and $X\subseteq X'$ be
sets of variables and $\B_RM_X\subseteq \B_RM_{X'}$ the
corresponding free bodies. If $u\in \B:=\B_RM_X$ and $v\in
X'\setminus X$, then $w:=u+v\in \B':= \B_RM_{X'}$ is free over $\B$,
i.e. there is a basis $X''$ of $\B'$ with $w\in X''\supseteq X$.
\end{fprop}

\Pf We will use Lemma \ref{charbasis} (c) to show that the set
$X'':= (X'\setminus \{v\})\cup \{w\}$ is a basis of $\B'$. First
note that $X''$ also generates $\B'$, thus $\B'=\B_RM_{X''}$.

Given $\va:X''\to H$ for a body $H$, we must extend this map to
$\bar{\va}:\B'\to H$. Let $\va': = \va\restr (X'\setminus \{v\})$
and note that the set $X'\setminus \{v\}= X''\setminus \{w\}$ is
independent. Thus if $\B_0:=\B_RM_{X'\setminus \{v\}}$, then $\va'$
extends to $\bar{\va'}: \B_0\to H$ by freeness, and from $u\in \B_0$
follows the existence of $u\bar{\va'}\in H$. We now define
$\bar{\va}$: if $\bar{\va}\restr \B_0 :=\bar{\va'}$, then
$\bar{\va}\restr (X'\setminus \{v\}) = \va'=\va\restr (X'\setminus
\{v\})$. Thus it remains to extend $\bar{\va'}$ to $\bar{\va}:\B'\to
H$ in such a way that $w\bar{\va}=w\va$. If $w\va=: h\in H$, then we
must have $h=w\bar{\va}= (u+v)\bar{\va}=u\bar{\va}+v\bar{\va}$.
Hence put $v\bar{\va} := h-u\bar{\va}=h-u\bar{\va'}$. Now
$\bar{\va}:\B'\to H$ exists, because $X'$ is free, $\va'\subseteq
\bar{\va}$ and $w\bar{\va}= (u+v)\bar{\va}=u\bar{\va}+ h -
u\bar{\va}= h=w\va$, thus $\va\subseteq \bar{\va}$ holds as
required.  \fine

\bigskip

The following corollaries (used several times for exchanging basis
elements) are immediate consequences of the last proposition.

\begin{corollary} \label{xbasis} Let $X$ be a basis for the body $\B$, $v\in X$
and $\B'$ be the subbody of $\B$ generated by $X\setminus \{v\}$ and
$w\in \B'$, then $X'=X\setminus \{v\}\cup \{v+w\}$ is another basis
for $\B$.
\end{corollary}

\begin{corollary} If $X\subseteq X'$ and $\B_RM_X \subseteq \B_RM_{X'}$,
then any basis of $\B_RM_X$ extends to a basis of $\B_RM_{X'}$.
\end{corollary}

\Pf If $X''$ is a basis of $\B_RM_X$, then it is left as an
exercise to show that $(X'\setminus X)\cup X''$ is a basis of
$\B_RM_{X'}$. \fine
\bigskip

The last corollaries have another implication.

\begin{corollary} \label{freechain} Suppose that $\B_\a \ (\a \le \d)$ is an
ascending, continuous chain of bodies such that $\B_{\a+1}$ is free
over $\B_\a$ for all $\a <\d$. Then $\B_\d$ is free over $\B_0$ and
if $\B_0$ is free, then $\B_\d$ is free as well.
\end{corollary}

The proof of the next lemma is also obvious. It follows by
application of the distributive law in $T^{bd}$ and collection of
summands with $p$.

\begin{lemma} Let $q(y)$ be a generalized polynomial
and $r\in R$. Then there is a polynomial $q'(y)$ such that
$q(y_1+ry_2) = q(y_1) + r q'(y_1,y_2)$.
\end{lemma}

\Proof By Lemma \ref{linearexpr} we can write $\s=\sum_{l<t}a_l\s_l$
with $(\s_l,\xbar)\in t(\tau^{sk},X)$ for any specified term
$(\s,\xbar)\in t(\tau^{bd},X)$. Thus it is enough to show
$q(y_1+ry_2) = q(y_1) + rq'(y_1,y_2)$ for generalized monomials $q$,
and this is obvious by iterated use of axiom (\ref{bodytheory})
(ii). \fine

\begin{lemma}\label{linearize} Let $X_0,X_1,X_2$ be pairwise disjoint infinite
sets, $\B_0:=\B(X_0)\subseteq \B:=\B(X_0\cup X_1\cup X_2)$ and
$q(y),q_1(y),q_2(y)$ polynomials over $\B_0$ such that
$$q(g+v_1+v_2) = q_1(v_1)+q_2(v_2)$$ for some $g\in \B_0, v_1\in
X_1,v_2\in X_2$. Then the following holds.
\begin{enumerate}
    \item $q(y)$ is a linear polynomial in $y$, i.e. $y$ appears
    at most once in every monomial.
    \item $q_1(y) - q_2(y)$ does not depend on $y$.
\end{enumerate}
\end{lemma}

\Pf (i) \  Write $q(y)=\sum_{i=1}^n m_i(y)$ as a sum of minimal
length of \index{generalized!monomial} generalized $R$--monomials
and suppose for contradiction that $y$ appears $n$ times in $m_1(y)$
with $n>1$. Also let, without loss of generality, $n$ be maximal for
the chosen monomial $m_1(y)$.

By the distributive law the monomials of the polynomial
$q(g+v_1+v_2)$ include those monomials induced by $m_1(y)$ replacing
all entries of the variable $y$ by arbitrary choices of $v_1$ and
$v_2$. Let $m'_1$ be one of these monomials. If there are further
such monomials $m'_i \ (i\le k)$ alike $m'_1$ arriving from this
substitution into monomials $m_i(y)$ of $q(y)$ with
$\sum_{i=1}^km'_i=0$, then replacing all $v_1$s and $v_2$s by $y$s
gives $\sum_{i=1}^km_i(y)=0$ contradicting the minimality of the
above sum. Thus $m'_1$ represents a true monomial (not canceled by
others) of $q(g+v_1+v_2)$, and as $n>1$ we may also assume that
$v_1$ and $v_2$ both appear in $m'_1$. This monomial does not exist
on the right--hand side of the equation in the lemma -- a
contradiction. Thus (a) holds.

(ii) \ First substitute in the given equation  $v_1:=y,v_2:=0$ and
$v_1:=0,v_2:=y$, respectively. Thus $q(g+y) = q_1(y) + c$ and
$q(g+y) = q_2(y) + c'$, where $c:=q_2(0), c' :=q_1(0) \in \B_0$.
Subtraction now yields $0=q_1(y) - q_2(y) + (c-c')$, thus $q_1(y) -
q_2(y)= c'-c$ does not depend on $y$, as required.  \fine

\bigskip

In order to establish the Step Lemmas, we next prepare some
preliminary results. Let $X_\o=\bigcup_{n\in\o} X_n$ be a strictly
increasing sequence of infinite sets $X_n$ of variables and fix a
sequence $v_n\in X_n\setminus X_{n-1}$ of elements ($n\in\o$).
Moreover, let $M_\a= M_{X_\a}$ be the skeleton and $\B_\a:=\B(X_\a)$
be the body generated by $X_\a$ for $\a \le \o$, respectively. Note
that by our identification $\B(X_\a)$ is an $R$-algebra and
restricting to the module structure $G_\a:={}_R\B(X_\a)$ is an
$R$-module, which is free by Lemma \ref{freeR}. Recall that $\bS=\{
p^n: n\in\o\}$ for some $p\in R$ (with $\bigcap_{n\in\o}p^nR=0$)
generates the $\bS$-topology on $R$-modules. Thus the $\bS$-topology
is Hausdorff on $G_\a$ and $G_\a$ is naturally an $\bS$-pure
$R$-submodule of its $\bS$-completion $\Ghat_\a$; we write
$G_\a\subseteq_*\Ghat_\a$ and pick particular elements $w_n\in
\Ghat_\o$. If $a_i\in \{0,1\}$ and $l_n\in \N$ is increasing, then
we define

\begin{eqnarray}\label{wn}w_n(v_n,l_n,a_n):=
w_n:= \sum\limits_{k\ge n}p^{l_k-l_n}a_kv_k \in \Ghat_\o
\end{eqnarray}

and easily check that
\begin{eqnarray}\label{wnn+1} w_n-p^{l_{n+1}-l_n}w_{n+1} = a_nv_n \in G_n \text{ for all }
n\in\o.
\end{eqnarray}

\begin{proposition} \label{prestep} Let $a_n\in \{0,1\}$ and $l_n\in \N$ be as above.
If $X_{\o+1} = X_\o \setminus \{v_n : a_n=1, n > 0\} \cup W$ with
$W=\{ w_n : n > 0\}$ and $\B_{\o+1} :=\B(X_{\o+1})$, $G_{\o+1}
={}_R\B_{\o+1}$, then the following holds
\begin{enumerate}
    \item  $G_\o\subseteq_*G_{\o+1} \subseteq_*\Ghat_\o$
    \item  $G_{\o+1}/G_\o$ is $p$-divisible, thus an $\bS^{-1}R$-module.
    \item  $X_{\o+1}$ is a basis of the (free) skeleton $M_{\o+1}=M_{X_{\o+1}}$.
    \item  The $R$-algebra $\B_{\o+1}$ is freely generated by the
    skeleton $M_{\o+1}$, thus $\B_{\o+1}= RM_{\o+1}$ and
    $G_{\o+1}=\bigoplus_{m\in M_{\o+1}} Rm$.
    \item  $\B(X_{\o+1})$ is free over $\B(X_n)$ (as body).
\end{enumerate}
\end{proposition}

\Pf Claim (i): \ Clearly $G_\o \subseteq_*\Ghat_\o$ and $G_{\o+1}
\subseteq\Ghat_\o$. From $w_n, w_{n+1}\in X_{\o+1}$, $a_n=1$ and
(\ref{wnn+1})  follows $v_n\in G_{\o+1}$. Hence $v_n\in G_{\o+1}$
for all $n\in\o$ and $G_\o\subseteq \B(X_{\o+1})=G_{\o+1}$ follows
at once. Thus $G_{\o+1}/G_\o\subseteq \Ghat_\o/G_\o$ and purity
($G_{\o+1} \subseteq_*\Ghat_\o$) follows if $G_{\o+1}/G_\o$ is
$p$-divisible. This is our next
\bigskip

Claim (ii): By definition of the body $\B(X_{\o+1})$, any element
$g\in G_{\o+1}$ is the sum of monomials in $X_{\o+1}$. If $w_n,
w_{n+1}$ are involved in such a monomial, then we apply
(\ref{wnn+1}) and get $w_n =pw_{n+1} + a_nv_n$, which is
$w_n\equiv p w_{n+1} \mod G_\o$. Let $m$ be the largest index of
$w_n$s which contributes to $g$. We can remove all $w_i$ of
smaller index $i < m$ and also write $w_m \equiv pw_{m+1}\mod
G_\o$. Thus $g+ G_\o$ is divisible by $p$ and $G_{\o+1}/G_\o$ is
an $\bS^{-1}R$-module.

\bigskip

Claim (iii): It is enough to show that $X_{\o+1}$ is free, because
$X_{\o+1}$ generates $M_{\o+1}$ by definition of the skeleton. First
we claim that $$X'=(X_\o \setminus \{v_n\})\cup \{w_n\} \text{ is
free.}$$ We apply the characterization of a basis by Lemma
\ref{charbasis} (ii). Let $(\s_1,\xbar),(\s_2,\xbar)\in
t^r(\tau^{bd},X)$ ($\xbar=\langle x_1,\dots,x_k\rangle$) be such
that
$$\s_1(y_1,\dots, y_k) =\s_2(y_1,\dots, y_k)$$  for some $y_i\in X'$ and
suppose that $y_1=w_n$ (there is nothing to show if $w_n$ does not
appear among the $y_i$s, because they are free; otherwise we relabel
the $y_i$s such that $y_1 =w_n$). Now we consider the above equation
as an element in $\Ghat_\o$ and note that the support
$[y_i]\subseteq X_\o$ of the elements $y_i, \ (i>1)$ is finite,
while $w_n$ has infinite support $\{ v_k : k>n \}\subseteq [w_n]$.
Thus we project $\s_i(y_1,\dots, y_k)$ onto a free summand from
$[w_n] \setminus \bigcup_{1<i\le k}[y_i]$ and $y_1$ can be replaced
by a free variable $v$ (over $y_2,\dots,y_k$) and $\s_1(v,y_2,\dots,
y_k) =\s_2(v,y_2,\dots, y_k)$ are the same. Hence the first claim
follows.

By the first claim and induction it follows that
\begin{eqnarray}\label{freeBn}
(X_\o\setminus \{v_1,\dots, v_n\})\cup \{w_1,\dots, w_n\} \text{
is free.}
\end{eqnarray}

Finally let $y_1,\dots, y_m$ be any finite subset of $X_{\o+1}$.
We may assume that $y_1,\dots, y_k\in W_n =\{w_i : i\le n\}$ and
$y_{k+1},\dots, y_m\in X_{\o+1}\setminus W$. Hence $$y_1,\dots,
y_m\in W_n\cup X_\o\setminus \{v_1,\dots,v_n\}$$ which is free by
the last claim (\ref{freeBn}), thus $X_{\o+1}$ is free and (iii)
follows.
\bigskip

(iv) is a consequence of (iii) and the definitions.

(v) We note that (by (iii)) the body $\B(X_{\o+1})$ is freely
generated by $X_{\o+1}$ and also (using (\ref{wnn+1})) by the
(free) set $X'= (X_{\o+1}\setminus\{w_1,\dots, w_n\})
\cup\{v_1,\dots, v_n\}$. However $X_n\subseteq X'$  which
generates $\B(X_n)$, hence (v) also follows. \fine

\bigskip

Throughout the remaining part of this section and Section \ref{S5}
we use the notations from Proposition \ref{prestep}. Moreover we
assume the following, where we view $\B_\o$ as an $R$-algebra.

\begin{eqnarray}\label{phi} \text{ Let } \va\in \End_RG_\o\setminus
\B_\o, \text{ with } G_n\va\subseteq G_n \text{ for all } n\in\o.
\end{eqnarray}

\begin{lemma}\label{catchw0} Let $\va$ be as in (\ref{phi}).
If $w_0\va\in G_{\o+1}$, then the following holds.
\begin{enumerate}
    \item There exist $m\in\o$ and a generalized polynomial $q_0(y)$
    over $\B_\o$ such that $w_0\va= q_0(w_m)$.
    \item There exists an $n_* > m$ such that $q_0$ is a polynomial over
    $\B_{n_*}$.
\end{enumerate}
\end{lemma}

\Pf (i)\ If $w_0\va\in G_{\o+1}$, then there exists some $m\in\o$
such that $w_0\va\in\B_RM_{X_{\o}\cup\{ w_0,\dots, w_{m}\} }$. Using
$w_i\equiv pw_{i+1} \mod G_\o$ it follows that
$$w_0\va\in\B_RM_{X_{\o}\cup\{ w_{m}\} },$$ and there is a
generalized polynomial $q_0(y)$ over $\B_\o$ such that $w_0\va=
q(w_m)$.

(ii)\ The coefficients of $q_0$ are in some $\B_{n_*}$ for some
$n_*>m$.
 \fine

\section{The three Step (or Stop) Lemmas}\label{S5}
We will use the notations from Proposition \ref{prestep} and
(\ref{phi}).

We begin with our first Step Lemma, which will stop $\va$ becoming
an endomorphism of our final module.

\begin{steplemma} \label{step1} Let $\va\in\End_RG$ be an endomorphism
as in (\ref{phi}) such that for all  $n\in\o \text{ there is }
g_n\in X_{n+1}\setminus G_n$  with  $g_n\va\notin \B_RM_{X_n\cup
\{g_n\}}$, and let $G_{\o+1}$ be defined with $w_n= w_n(g_n,l_n,1)$
as in (\ref{wn}) for suitable elements $l_n\in\o$. Then $\va$ does
not extend to an endomorphism in $\End_RG_{\o+1}$.
\end{steplemma}

\Pf We define inductively an ascending sequence $l_n\in \o$:\\
If $C\subseteq G_\o$ is a submodule, then the $p$-closure of $C$
is defined by $\Cbar= \bigcap_{n\in\o}(p^nG_\o + C)$. It is the
closure of $C$ in the $\bS$-adic topology, which is Hausdorff on
$G_\o$ (i.e. $\bigcap_{n\in\o}p^nG_\o=0$). In particular $C=\Cbar$
if $C$ is a summand of $G_\o$ (e.g. $C=0$ is closed).

By hypothesis we have $g_n\va\notin \B_RM_{X_n\cup \{g_n\}}$ and
$\B_RM_{X_n\cup \{g_n\}}$ is a summand of $G_\o$, so it is closed in
the $\bS$-topology. There is an $l\in\o$ such that $g_n\va\notin
\B_RM_{X_n\cup \{g_n\}}+ p^lG_\o$. If $l_{n-1}$ is given, we may
choose $l=l_n$ such that $l_n> 3l_{n-1}$. We will ensure (just
below) that $G_\o\subseteq_* G_{\o+1}$ is $\bS$-pure, thus
$p^{l_n}G_{\o+1}\cap G_\o\subseteq p^{l_n}G_\o$, and it follows:
\begin{eqnarray}\label{forcontra} \text{ There is a sequence } l_n\in \o \text{ with }
l_{n+1}>3l_n \text{ and } g_n\va\notin \B_RM_{X_n\cup \{g_n\}}+
p^{l_n}G_{\o+1}.\end{eqnarray}

Hence $G_{\o+1}$ is well defined and Proposition \ref{prestep}
holds; in particular $G_\o\subseteq_* G_{\o+1}$ and $G_\o$ is dense
in $G_{\o+1}$ ($\bar G_\o =G_{\o+1}$). Suppose for contradiction
that $\va\in\End_RG_\o$ extends to an endomorphism of $G_{\o+1}$;
this extension is unique, and we call it also
$\va\in\End_RG_{\o+1}$. In particular $w_0\va\in G_{\o+1}$; by Lemma
\ref{catchw0}(i)(ii) there is a polynomial $q_0(y)$ with
coefficients in $\B_{n_*}$ for some $n_*\in\o$ and with
$w_0\va=q_0(w_m)$ for some $m\in\o$. We choose $n>\max\{n_*,m\}$ and
use (\ref{wn}) to compute $w_0$:
$$w_0= \sum\limits_{i=1}^n p^{l_i-l_0}g_i + p^{l_{n+1}-l_0}w_{n+1}.$$

Application of $\va$ gives
$$w_0\va\equiv \sum\limits_{i=1}^{n-1} p^{l_i-l_0}(g_i\va) +  p^{l_n-l_0}(g_n\va) \mod p^{l_{n+1}-l_0}G_{\o+1}.$$
If $i<n$, then $g_i\in G_n$ and $g_i\va\in G_n$ by the choice of
$\va$. The last equality becomes $w_0\va\equiv p^{l_n-l_0}(g_n\va)
\mod p^{l_{n+1}-l_0}G_{\o+1}+G_n$, hence
$$q_0(w_0) \equiv p^{l_n-l_0}(g_n\va) \mod (p^{l_{n+1}-l_0}G_{\o+1}+G_n).$$

Finally we determine $p^{l_n-l_0}(g_n\va)$ in terms of
$\B_RM_{X_n\cup \{g_n\}}$:

\relax From $w_0\va=q_0(w_m)$, $n>m$ and the definition of $w_m$ in
(\ref{wnn+1}) we get $$w_m = \sum\limits_{i=m}^{n-1} p^{l_i-l_m}g_i
+ p^{l_n -l_m}g_n + p^{l_{n+1}-l_m}w_{n+1},$$ thus
$$q_0(w_m) \equiv q_0(\sum\limits_{i=m}^{n-1} p^{l_i-l_m}g_i
+ p^{l_n -l_m}g_n) \mod p^{l_{n+1}-l_m}G_{\o+1},$$ and
$$
p^{l_n-l_0}(g_n\va) \equiv q_0(w_m)  \equiv
q_0(\sum\limits_{i=m}^{n-1} p^{l_i-l_m}g_i + p^{l_n -l_m}g_n) \mod
(p^{l_{n+1}-l_m}G_{\o+1}+ G_n).
$$

Now we use (again) that $g_i\in G_n$ for all $i<n$. The last
equation reduces to $ p^{l_n-l_0}(g_n\va) \in  \B_RM_{X_n\cup
\{g_n\}} + p^{l_{n+1}-l_m}G_{\o+1}$, hence $g_n\va \in  \B(G_n,g_n)
+ p^{l_{n+1}-l_m-l_n}G_{\o+1}$. Note that $l_{n+1} > 3l_n$ by the
choice of the $l_n$s, hence $l_{n+1}-l_m-l_n>l_n$, so we get a
formula
$$g_n\va \in \B_RM_{X_n\cup \{g_n\}} + p^{l_n}G_{\o+1}$$ that contradicts
(\ref{forcontra}) and the Step Lemma \ref{step1} follows. \fine

\begin{steplemma} \label{step2} Let $\va\in\End_RG$ be an endomorphism
as in (\ref{phi}). Moreover suppose there are elements $u_n,g_n\in
X_{n+1}\setminus G_n$ (for each $n\in \o$) with $u_n\va= q^1_n(u_n)$
and $g_n\va= q^2_n(g_n)$, where $q^1_n,q^2_n$ are polynomials over
$\B_0$ such that
$$q_n^1(y)- q_n^2(y)\notin \B_0 \text{ i.e. $y$ appears in the difference.} $$
If $G_{\o+1}$ is defined with $w_n= w_n(g_n+u_n,l_n,1)$ as in
(\ref{wn}) for suitable elements $l_n\in\o$, then $\va$ does not
extend to an endomorphism in $\End_RG_{\o+1}$.
\end{steplemma}

\Pf Let $c_k:= g_k+u_k$. The set $X':=(X_\o \setminus \{u_k | \ k<
\o\}) \cup \{c_k |\ k<\o\}$ is now a basis of $\B_\o$ by Corollary
\ref{xbasis}, thus Proposition \ref{prestep} applies and $G_{\o+1}$
is well--defined. By definition of $w_n$ and (\ref{wnn+1}) we have
$p^{l_{n+1}-l_n}w_{n+1}+c_n=w_n$, and as in the proof of Step Lemma
\ref{step1} we get
\begin{align*}
w_0\va=q_0(w_m) \Longrightarrow \sum\limits_{k\ge 0}
p^{l_k-l_0}(c_k\va)= q_0(\sum\limits_{k\ge m} p^{l_k-l_m}c_k)
\end{align*}
 for $n_*,m,q_0(y)$ as in Step Lemma
 \ref{step1}. Furthermore,
\begin{align*}
c_k\va = (g_k+u_k)\va = g_k\va + u_k\va = q_k^1(g_k)+q^2_k(u_k).
\end{align*} Thus $$\sum\limits_{k\ge 0}
p^{l_k-l_0}(q_k^1(g_k)+q^2_k(u_k))= q_0(\sum\limits_{k\ge m}
p^{l_k-l_m}(g_k+u_k)),$$ where $q_k^1(g_k) \in {_R\B}_{X_0
\cup\{g_k\}}$ and $q_k^2(u_k) \in {_R\B}_{X_0 \cup\{u_k\}}$, and
arguments similar to Lemma \ref{linearize} apply: in every monomial
of $q_0(y)$ the variable $y$ appears at most once as there are no
mixed monomials on the left--hand side, and the same holds for
$q^1_k(y), q^2_k(y)$. Furthermore, the variable $y$ does not appear
in $q^1_k(y)-q^2_k(y)$, which contradicts our assumption on the
$q^i_k$s.  \fine

The next lemma is the only place where we will use that $R$ is
$\Sigma\bS$-incomplete in order to find a sequence $a_n\in \{0,1\}$,
see Definition \ref{sigmaS}. Recall that this condition follows by
Corollary \ref{corGM} if the $\bS$-ring is a direct sum of
$\bS$-invariant subgroups of size $<\Cont$. Hence it will be
sufficient if $R^+$ is free and $\bS$ defines the usual $p$-adic
topology on $R$.

\begin{steplemma} \label{step3} Let $R$ be a $\Sigma\bS$-incomplete
$\bS$-ring, let $\va\in\End_RG$ be an endomorphism as in (\ref{phi})
and let $q=q(y)$ be a polynomial in $y$ with coefficients in $\B_0$
such that $g\va - q(g)\in G_0$ for all $g\in G$. Moreover suppose
that for all $n\in \o$ there are elements
$$g_n\in X_{n+1}\setminus G_n \text{ such that } g_n\va - q(g_n)\ne
0.$$ If $G_{\o+1}$ is defined with $w_n= w_n(g_n\va-q(g_n),l_n,1)$
as in (\ref{wn}) for suitable elements $l_n\in\o$, then $\va$ does
not extend to an endomorphism in $\End_RG_{\o+1}$.
\end{steplemma}

\Pf Choose $g_n\in X_{n+1}$ as in the Lemma, and put
$h_n=g_n\va-q(g_n)\ne 0$. By assumption on $\va$ and $q$ it follows
that $h_n\in G_0$. Let $w_n = w_n(g_n,n,a_n)$ be defined as in
(\ref{wn}) for a suitable sequence of elements $a_n\in\{0,1\}$ and
$l_n=n$ for all $n\in\o$. We define again $G_{\o+1}$ as in
Proposition \ref{prestep} using the new choice of elements $w_n$.
Note that $G_0$ is a free $R$-module. By the assumption that $R$ is
$\Sigma\bS$-incomplete there is a sequence $a_n\in\{0,1\}$ with
$\sum_{k\in\o}p^ka_kh_k\notin G_0$. However
$\sum_{k\in\o}p^ka_kh_k\in \widehat G_0$ by the choice of $h_n$,
hence $\sum_{k\in\o}p^ka_kh_k\notin G_{\o+1}$ by definition of
$G_{\o+1}$. Recall $w_0= \sum_{k\in\o}p^ka_kg_k$ and suppose that
$w_0\va\in G_{\o+1}$. We compute
$$w_0\va=(\sum_{k\in\o}p^ka_kg_k)\va = \sum_{k\in\o}p^ka_k(g_k\va)$$
and $$\sum_{k\in\o}p^ka_kh_k= \sum_{k\in\o}p^ka_k(g_k\va-q(g_k))=
\sum_{k\in\o}p^ka_kg_k - \sum_{k\in\o}q(g_k)p^ka_k=w_0-
\sum_{k\in\o}p^ka_kq(g_k).$$ From $\sum_{k\in\o}p^ka_kh_k \notin
G_{\o+1}$ follows $\sum_{k\in\o}p^ka_kq(g_k)\notin G_{\o+1}$.
However, by the definition of bodies $\B_\a$,  the map taking
$g\arr q(g)$ for any $g\in G_{\o+1}$ is an endomorphism of
$G_{\o+1}$, and also $w_0=\sum_{k\in\o}p^ka_kg_k\in G_{\o+1}$,
hence $\sum_{k\in\o}p^ka_kq(g_k)\in G_{\o+1}$ is a contradiction.
We deduce $w_0\va\notin G_{\o+1}$ and $\va$ does not extend to an
endomorphism of $G_{\o+1}$. \fine

\section{ Constructing Generalized $\boldsymbol{E(R)}$-Algebras}

\begin{lemma}\label{freeele} Let $\k$ be a regular, uncountable cardinal and
$\B=\bigcup_{\a\in\k} \B_\a$ a $\k$-filtration of bodies. Also let
$G_\a={}_R\B_\a$ and $G={}_R\B$. Then the following holds for any
$\va\in \End_RG$.
\begin{enumerate}
 \item If there is $g\in G$ such that $g\va \notin (\B_\a)_{\{g\}}$, then
 there is also $h\in G$ free over $\B_\a$ such that $h\va \notin
 (\B_\a)_{\{h\}}$.
 \item If there are $g\in G$ and a polynomial $q(y)$ over $\B_\a$
 such that $g\va - q(g)\notin (\B_\a)_{\{g\}}$, then there is also
 $h\in G$ free over $\B_\a$ such that $h\va -
 q(h)\notin (\B_\a)_{\{h\}}$.
\end{enumerate}
\end{lemma}

\Pf If $g\in G$ satisfies the requirements in (i) or (ii),
respectively, then choose any element $h'\in G$ which is free over
$\B_\a$. If $h'$ also satisfies the conclusion of the lemma, then
let $h=h'$ and the proof is finished. Otherwise let $h=h'+g$ which
is  also free over $\B_\a$ by Proposition \ref{fproper}. In this
case $h'\va\in (\B_\a)_{\{h'\}}$ or $h'\va-q(h')\in
(\B_\a)_{\{h'\}}$, respectively. It follows $h\va\notin
(\B_\a)_{\{h\}}$ or $h\va-q(h)\notin (\B_\a)_{\{h\}}$, respectively.
\fine

\bigskip

The next lemma is based on results of the last section concerning
the Step Lemmas and Lemma \ref{freeele}. We will construct first the
$\k$-filtration of $\B_\a$s for application using $\Diam_\k E$ for
some non-reflecting subset $E\subseteq \k^o$. Recall that $\Diam_\k
E$ holds for all regular, uncountable, not weakly compact cardinals
$\k$ and non-reflecting subsets $E$ in $V=L$.

\bigskip
\noindent
 {\sl Construction of a $\k$-filtration of free bodies:}
Let $\{\va_\rho: \rho\in E\}$ be the family of Jensen functions
given by $\Diam_\k E$. The body $\B$ and the $R$-module ${}_R\B$
will be constructed as a $\k$--filtration $\B =\bigcup\limits_{\a
\in \k} \B_\a$ of bodies. We choose
$$|{\B_\a}| =|\a|+|R|=|\B_{\a +1} \setminus \B_\a|$$
and fix for each $\a\in E$ a strictly increasing sequence
$$\a_n\in \a\setminus E \text{  with } \sup_{n\in\o}\a_n = \a.$$
This is possible, because $E$ consists of limit ordinals cofinal
to $\o$ only and we can pick $\a_n$ as a successor ordinal. We
will use the same Greek letter for a converging sequence and its
limit, so the elements of the sequence only differ by the suffix.

As $E$ is non-reflecting, we also may choose a strictly
increasing, continuous sequence $\a_\nu, \ \nu \in \cf(\a)$ with
$$\sup_{\nu \in \cf(\a)}\a_\nu =\a \text{ and } \a_\nu\in \a\setminus E$$
if $\cf(\a) > \o$. This is crucial, because the body $\B_\a$ of
the (continuous) $\k$--filtration of $\B$ must be free in order to
proceed by a transfinite construction. This case does not occur
for $\k =\ale $.

Using Lemma \ref{step1},  Lemma \ref{step2} and Lemma \ref{step3}
inductively, we define the body structure on $\B_\nu$. We begin
with $\B_0 = 0$, and by continuity of the ascending chain the
construction reduces to an inductive step passing from $\B_\nu$ to
$\B_{\nu +1}$. We will carry on our induction hypothesis of the
filtration at each step. In particular the following three
conditions must hold.
\begin{enumerate}
 \item $\B_\nu$  is a free body.
 \item If $\rho \in \nu \setminus E$, then $\B_\nu$ is a free body over $\B_\rho$.
 \item If $\rho \in E$ then let $\rho_n \ (n\in\o)$ be the given
 sequence with $\sup_{n\in \o} \rho_n = \rho$. Suppose that the
 hypothesis of one of the three step lemmas, Lemma \ref{step1} or
 Lemma \ref{step2} or Lemma \ref{step3} holds for $G_n =
{}_R\B_{\rho_n}$ ($n\in\o$). We identify $\B_{\rho+1}$ with
$\B_{\o+1}$ from the Step Lemmas (so $\va_\rho$ does not extend to
an endomorphism of $\B_{\rho +1}$).
\end{enumerate}

Following these rules we step to $\nu +1$: If the hypotheses of
condition (iii) are violated, for instance, if $\nu \notin E$ we
choose $\B_{\nu +1} := (\B_\a)_{\{ v_\a\}}$ adding any new free
variable $v_\nu$ to the body. However, next we must check that these
conditions (i) to (iii) can be carried over to $\nu +1$. If the
hypotheses of condition (iii) are violated, this is obvious. In the
other case the Step Lemmas are designed to guarantee:

Condition (i) is the freeness of $\B_{\o+1}$ in Proposition
\ref{prestep}. Condition (ii) needs that $\B_{\nu+1}$ is a free
body over $\B_\rho$. However $\B_\rho\subseteq \B_{\nu_n}$ for a
large enough $n\in \o$. Hence (ii) follows from freeness of
$\B_{\nu_n}$ over $\B_\rho$ (inductively) and $\B_{\nu+1}$ over
$\B_{\nu_n}$ (by Proposition~ \ref{prestep} and Corollary
\ref{freechain}).

In the limit case $\g$ we have two possibilities:  If $\cf(\g) = \o$
then $\sup_{n\in \o} \g_n = \g$, hence $\B_\g
=\bigcup\limits_{n\in\o}\B_{\g_n}$ and $\B_\g$ is a free body with
the help of (i) and (ii) by induction, see Corollary
\ref{freechain}.  If $\cf(\g) > \o$, then by our set theoretic
assumption ($E$ is non--reflecting) we have a limit $\sup_{\a\in
\cf(\g)}\g_\a = \g$ of ordinals not in $E$. The union of the chain
$\B_\g = \bigcup\limits_{\a\in \cf(\g)}\B_{\g_\a}$ by (i) and (ii)
is again a free body, see Corollary \ref{freechain}. Thus we proceed
and obtain $\B=\bigcup\limits_{\nu\in\k} \B_\nu$ which is a
$\k$--filtration of free bodies. It remains to show the

\begin{mlemma} \label{mainl} Assume $\Diam_\k E$. Let $\k$ be a regular,
uncountable cardinal and $\B=\bigcup_{\a\in\k} \B_\a$ be the
$\k$-filtration of bodies just constructed. Also let $G_\a={}_R\B_\a
$ and $G={}_R\B$. Suppose that $\va\in \End_RG$ does not satisfy the
following conditions (i) or (ii) for any $ \a\in\k$ and any
polynomial $q(y)$ over $\B_\a$:
\begin{enumerate}
\item There is $g\in G$ such that $g\va \notin {}_R(\B_{\b})_{\{ g\}}$. \item
There is $g\in G$ such that $g\va - q(g)\notin {}_R(\B_{\b})_{\{
g\}}$
\end{enumerate}
Then $\va$ is represented in $\B$.
\end{mlemma}

\Pf Suppose for contradiction that $\va$ is not represented in
$\B$. Let $E\subseteq \k^o$ be given from $\Diam_\k E$, let
$\{\va_\d: \d\in E\}$ be the family of Jensen functions and define
a stationary subset $E'_\va=\{\d \in E: \va\restr G_\d=\va_\d\}$.
Note that $C=\{\d\in\k: G_\d\va\subseteq G_\d\}$ is a cub, thus
$E_\va:=E'_\va\cap C$ is also stationary.

As a consequence we see that there is $\d\in E_\va$ satisfying one
of the following conditions.
\begin{enumerate}
 \item For every $\a <\d\in E_\va$ there is $g\in G_\d$ such that $g\va
 \notin {}_R(\B_{\a})_{\{ g\}}$.
 \item There is $\a < \d\in E_\va$ such that $g\va\in  {}_R(\B_{\a})_{\{ g\}}$ for
 all $g\in G_\d$ (not case (i)) but for every $\a < \d$ and every
 polynomial $q(y)$ over $\B_\a$ represented by an endomorphism of $G$
 there is $g\in G_\d$ with $g\va-q(g)\notin G_\a$.
 \item There is $\a < \d\in E_\va$ such that $g\va\in  {}_R(\B_{\a})_{\{ g\}}$ and there
 is a polynomial $q(y)$ over $\B_\a$ with $g\va-q(g)\in G_\a$ for
 all $g\in G_\d$ (neither (i) nor (ii) holds), but $\va$ is not
 represented by $\B$. Thus there are a sequence $\d_n<\d$ ($n\in \o$)
 with $\sup_{n\in\o}\d_n=\d$ and $g_n\in G_{\d_n} $ such that
 $g_n\va-q(g_n)\ne 0$ for all $n\in\o$.
\end{enumerate}

By Lemma \ref{freeele} we may assume that the elements $g$
existing by (i) and (ii), respectively, are free over $\B_\a$.
Moreover, if $g\in G$, then by $\cf\k >\o$ we can choose $\d\in
E_\va$ such that $g\in G_\d$.

If (i) holds, then we can choose a proper ascending sequence
$\d_n\in E_\va$ with $\sup_{n\in\o} \d_n=\d$ and elements $g_n\in
G_{\d_{n+1}}$ such that $g_n$ is free over $\B_{\d_n}$ and
 $$ g_n\va\notin {}_R(\B_{\d_n})_{\{ g_n\}} \text{ for all } n\in\o. $$

We identify $G_{\d_n}$ with $G_n$ in Step Lemma \ref{step1} and
note that (by $\d_n\in E_\va$) $\va\restr G_n$ is an endomorphism
with $G_n\va\subseteq G_n$ which is predicted as a Jensen
function. By construction of $G_{\d+1}$ (as a copy of $G_{\o+1}$
from Lemma \ref{step1}), the endomorphism $\va\restr G_\d$ does
not extend to $\End_RG_{\d+1}$. However $\va\in \End G$, thus
$G_{\d+1}\va\subseteq G_\a$ for some $\a<\k$. Finally note that
$G_{\d+1}$ is the $\bS$-adic closure of $G_\d$ in $G$ because
$G_\d$ is $\bS$-dense in $G_{\d+1}$ and $G_{\d+1}$ is a summand of
$G_\a$ hence $\bS$-closed in $G_\a$. We derive the contradiction
that indeed $\va\restr G_{\d+1}\in \End_RG_{\d+1}$. Hence case (i)
is discarded.

Now we turn to case (ii). Suppose that (ii) holds (so condition (i)
is not satisfied). In this case there is an ascending sequence
$\d_n\in E_\va$ with $\sup_{n\in\o} \d_n=\d$ as above and there are
free elements $g_n,u_n\in G_{\d_{n+1}}$ (also free over $\B_{\d_n}$
and polynomials $q_n^1,q_n^2$ over $\B_{\d_0}$ such that $u_n\va=
q^1_n(u_n)\ne g_n\va= q^2_n(g_n)$. Moreover, the polynomials
$q^1_n(y) - q^2_n(y)$ are not constant over $\B_{\d_0}$. Step Lemma
\ref{step2} applies and we get a contradiction as in case (i). Thus
also case (ii) is discarded.

Finally suppose for contradiction that (iii) holds (so (i) and (ii)
are not satisfied). There are $\a <\k$ and $q(y)$ a polynomial over
$\B_\a$ such that $g\va-q(g)\in G_\a$ for all $g\in G$. The
polynomial $q(y)$ is represented by an endomorphism of $G_\a$.
Moreover (from (iii)) we find $g_n\in G_{\d_{n+1}}$ free over
$G_{\d_n}$ for a suitable sequence $\d_n$ with $\sup_{n\in\o}
\d_n=\d$ such that $g_n\va - q(g_n)\ne 0$. We now apply Step Lemma
\ref{step3}; the argument from case (i) gives a final contradiction.
Thus the Main Lemma holds. \fine

{\sl Proof of Main Theorem \ref{mainth}.} Let $\B$ be the body
over the $\bS$-ring $R$ constructed at the beginning of this
section using $E$ as in Theorem \ref{mainth}; moreover let
$G={}_R\B$. Thus $\size{\B}=\k$ and by the construction and
Proposition \ref{prestep} any subset of size $<\k$ is contained in
an $R$-monoid-algebra of cardinality $<\k$; the algebra $\B$ is
the union of a $\k$-filtration of free bodies $\B_\a$. By Lemma
\ref{mainl} every element $\va\in \End_RG$
 is represented (by a polynomial $q(y)$ with
coefficients from $\B$); see Definition \ref{reppol}. Thus $g\va
=q(g)$ for all $g\in G$ and $\va=q(y)\in \B$. It follows that
$\B=\End_RG$ is the $R$-endomorphism algebra of $G$.

Finally recall that there is $\B_\a\subseteq \B$ which is an
$R$-monoid-algebra over a non-commutative monoid from Observation
\ref{freemono}. Thus $\B$ cannot be commutative either and Theorem
\ref{mainth} is shown. \fine

\goodbreak

\noindent R\"udiger G\"obel \\ Fachbereich Mathematik\\
  Universit\"at Duisburg-Essen\\ D 45117 Essen, Germany \\
  {\small e-mail: r.goebel@uni-essen.de}

 \noindent Saharon Shelah\\
Institute of Mathematics,
\\Hebrew University, Jerusalem, Israel \\
and Rutgers University, New Brunswick, NJ, U.S.A \\
{\small e-mail: shelah@math.huji.ac.il}

\end{document}